\documentclass[a4paper]{article}
\usepackage{amsthm}
\usepackage{amsmath}
\usepackage[margin=1.2in]{geometry}
\usepackage{amsfonts}
\usepackage{amssymb}
\usepackage{latexsym}
\usepackage{graphicx}
\usepackage{color}

\usepackage{tabularx}
\usepackage{enumitem}
\usepackage{ae}
\usepackage[T1]{fontenc}

\newtheorem{Theorem}{Theorem}[section]
\newtheorem{Proposition}[Theorem]{Proposition}
\newtheorem{Construction}[Theorem]{Construction}
\newtheorem{Corollary}[Theorem]{Corollary}

\newtheorem{Remark}[Theorem]{Remark}
\newtheorem{Lemma}[Theorem]{Lemma}
\newtheorem{Example}[Theorem]{Example}

\date{}

\begin{document}

\begin{center}
{\Large\bf Partitionable sets, almost partitionable sets and their applications \footnote{Partially supported by Italian Ministry of Education under Grant PRIN 2015 D72F16000790001 (S. Costa). Supported by NSFC under Grant 11971053 (Y. Chang), NSFC under Grant 11871095 (T. Feng), NSFC under Grant 11771227, and the Fundamental Research Funds for the Provincial Universities of Zhejiang under Grant SJLY2020008 (X. Wang).}}

\vskip12pt

Yanxun Chang$^a$, Simone Costa$^b$, Tao Feng$^a$, Xiaomiao Wang$^c$\\[2ex]
{\footnotesize $^a$Department of Mathematics, Beijing Jiaotong University, Beijing 100044, P. R. China}\\
{\footnotesize $^b$Dipartimento DICATAM, Universit\`a degli Studi di Brescia, Via Branze 43 - 25123 Brescia, Italy}\\
{\footnotesize $^c$School of Mathematics and Statistics, Ningbo University, Ningbo, Zhejiang 315211, P. R. China}\\
{\footnotesize yxchang@bjtu.edu.cn, simone.costa@unibs.it, tfeng@bjtu.edu.cn, wangxiaomiao@nbu.edu.cn}
\vskip12pt

\end{center}

\begin{abstract}
This paper introduces almost partitionable sets to generalize the known concept of partitionable sets. These notions provide a unified frame to construct $\mathbb{Z}$-cyclic patterned starter whist tournaments and cyclic balanced sampling plans excluding contiguous units. The existences of partitionable sets and almost partitionable sets are investigated. As an application, a large number of optical orthogonal codes achieving the Johnson bound or the Johnson bound minus one are constructed.
\end{abstract}



\section{Introduction}

The concept of partitionable sets was introduced in \cite{ZC2009} to construct cyclic balanced sampling plans excluding contiguous units (briefly BSEC) (see Section 2.3 for more details). Throughout this paper, we assume that $(G,+)$ is a group in which the identity is denoted by $0$. Sets and multisets will be denoted by curly braces $\{\ \}$ and square brackets $[\ ]$, respectively.

A partitionable set is a set of ordered pairs from a group such that additional properties are satisfied. More formally, let $(G,+)$ be a group of order $v\equiv 1\pmod{4}$. A {\em partitionable set in $G$, briefly PS$(G)$, is a set ${\cal S}$ of $(v-1)/4$ ordered pairs from $G$ such that
\begin{enumerate}
\item[$(1)$] $\bigcup\limits_{(x,y)\in {\cal S}} \pm\{x,y\}=G\setminus \{0\};$
\item[$(2)$] $\bigcup\limits_{(x,y)\in {\cal S}} \pm\{x-y,x+y\}=G\setminus \{0\}$.
\end{enumerate}
When $G$ is an abelian group, every pair in $\cal S$ can be seen as an unordered one and so we usually write $\{x,y\}$ instead of $(x,y)$.} If $G=\mathbb{Z}_v$, the additive group of integers modulo $v$, then the notation of PS$(v)$ is used.

The definition of partitionable sets requires $v\equiv 1 \pmod{4}$. It is natural to ask whether a similar definition can be given for $v\equiv 3 \pmod{4}$. Given a group $(G,+)$ of order $v\equiv 3 \pmod{4}$ and two nonzero elements $\alpha$ and $\beta$ of $G$, a $(G,\alpha,\beta)$-{\em almost partitionable set}, briefly APS$(G,\alpha,\beta)$, is a set $\cal S$ of $(v-3)/4$ ordered pairs from $G$ such that
\begin{enumerate}
\item[$(1)$] $\bigcup\limits_{(x,y)\in {\cal S}} \pm\{x,y\}=G\setminus\{0,\alpha,-\alpha\};$
\item[$(2)$] $\bigcup\limits_{(x,y)\in {\cal S}} \pm\{x-y,x+y\}=G\setminus\{0,\beta,-\beta\}$.
\end{enumerate}
When $G$ is an abelian group, we write $\{x,y\}$ instead of $(x,y)$. If $G=\mathbb{Z}_v$, we will use the notation of APS$(v,\alpha,\beta)$ or, if it is not important to highlight the roles of $\alpha$ and $\beta$, of APS$(v)$.

\begin{Example}\label{eg:APS}
An APS$(27,3,6)$: $\{ 1, 4 \}$, $\{ 2, 12 \}$, $\{ 5, 13 \}$, $\{ 6, 10 \}$, $\{ 7, 8 \}$, $\{ 9, 11 \}$.
\end{Example}

We observe that (see Proposition \ref{PSandZCPS}) when $v\equiv 1 \pmod{4}$, a PS$(v)$ exists if and only if there exists a $\mathbb{Z}$-cyclic patterned starter whist tournament of order $v$; when $v\equiv 3 \pmod{4}$, an APS$(v,\alpha,\alpha)$ for some $\alpha\in \mathbb Z_v$ exists if and only if there exists a $\mathbb{Z}$-cyclic patterned starter whist tournament of order $v+1$. The whist tournament problem has a history of more than one hundred years and much work has been done by many authors, but the existence of $\mathbb{Z}$-cyclic whist tournaments is still far from being settled (see Section 2.1). On one hand, Zhang and Chang \cite{ZC2009} initialized the study of PSs in 2009 and showed that any PS$(v)$ with $v\equiv 9\pmod {12}$ and $v\not\equiv 81\pmod{108}$ cannot exist. They also conjectured that no PS$(v)$ with $v\equiv 81\pmod{108}$ exist. On the other hand, Abel, Anderson and Finizio \cite{AAF2011} in 2011 made the same conclusion and the same conjecture for $\mathbb{Z}$-cyclic patterned starter whist tournament of order $v$. We have noticed the equivalence between these concepts only very recently and this is one of the reasons we are interested in PSs and APSs.

Another motivation comes from an application of PSs and APSs to optical orthogonal codes. We note that Buratti's silver ratio construction for optical orthogonal codes in \cite{b18} implies an  APS$(p,\alpha,\beta)$ with a prime $p\equiv 7 \pmod{8}$ (see Lemma \ref{Marco}). This observation makes us believe that the construction of PSs and APSs is a subject worth of deep investigation.

The rest part of this paper is organized as follows. Section 2.1 is devoted to exploring connections among PSs, APSs and various kinds of whist tournaments. Based on the relation between $\mathbb{Z}$-cyclic directed whist tournaments and cyclic difference matrices, Section 2.2 gives the fact that a PS$(v)$ implies a cyclic difference matrix of order $v$ with five rows. The necessary conditions for the existence of almost partitionable sets are discussed in Section 3. It will be shown that an APS$(v,\alpha,\beta)$ exists only if (see Lemma \ref{nec})
\begin{eqnarray}\label{eq:nec}
2\alpha^2-\beta^2\equiv
\left\{
\begin{array}{ll}
    \frac{v}{3} \ ({\rm mod }\ v), & \hbox{ \rm{if} $v\equiv 3\ ({\rm mod }\ 12)$,}\\[0.4em]
	0\ ({\rm mod }\ v), & \hbox{ \rm{if} $v\equiv 7,11 \ ({\rm mod }\ 12)$.}
\end{array}
\right.
\end{eqnarray}
Section 4 presents a direct construction for APS$(p,\alpha,\beta)$ with a prime $p\equiv 7 \pmod{8}$ and Section 5 gives recursive constructions for APSs by introducing partial partitionable sets. We will show that if there exist a PS$(v)$ and an APS$(u,\alpha,\beta)$ with $u\equiv 7,11\pmod{12}$, then there exists an APS$(vu,v\alpha,v\beta)$ (see Corollary \ref{cor:recursive}). As a corollary, infinite families of APS$(v,\alpha,\beta)$s are produced for any $\alpha$ and $\beta$ satisfying the basic necessary condition (see Corollary \ref{cor:APS from recur}).

In Section 6, we give a ``Kramer-Mesner type'' construction for PSs and APSs. It will be shown (see Theorem \ref{thm:<300}) that an APS$(v,\alpha,\beta)$ exists for any $v\equiv 3\pmod{4}$, $v<300$ and $\alpha$, $\beta$ satisfy the necessary condition \eqref{eq:nec}. Thus it can be conjectured that the necessary conditions for the existence of an APS$(v,\alpha,\beta)$ shown in Lemma \ref{nec} are also sufficient. On the other hand, our results on PSs and APSs allow us to obtain new $\mathbb{Z}$-cyclic patterned starter whist tournaments. Finally, as applications of partitionable sets and almost partitionable sets, we construct a large number of optical orthogonal codes achieving the Johnson bound or the Johnson bound minus one in Section 7. A conclusion is given in Section 8.

\section{Interplay with various designs}

\subsection{$\mathbb{Z}$-cyclic patterned starter whist tournaments}

In this section we will provide a link between the concepts of (almost) partitionable sets and whist tournaments.



A {\em whist tournament} on $v$ players, briefly Wh$(v)$, for $v=4n$ $($or $4n+1)$, is a schedule of games $(a, b, c, d)$ where the unordered pairs $\{a, c\}$, $\{b, d\}$ are called {\em partners} and the pairs $\{a, b\}$, $\{c, d\}$, $\{a, d\}$, $\{b, c\}$ are called {\em opponents}, such that
\begin{enumerate}
\item[$(1)$] the games are arranged into $4n-1$ $($or $4n+1)$ rounds,
each of $n$ games;
\item[$(2)$] each player plays in exactly one game in each round (or all rounds but one);
\item[$(3)$] each player partners every other player exactly once;
\item[$(4)$] each player has every other player as an opponent exactly twice.
\end{enumerate}

A Wh$(4n)$ is a particular type of $(4n,4,3)$ resolvable balanced incomplete block designs, and a Wh$(4n+1)$ is a particular type of $(4n+1,4,3)$ near resolvable balanced incomplete block designs. The whist tournament problem was introduced by Moore \cite{Moore} in 1896. Since 1970's, it has been known that a Wh$(v)$ exists if and only if $v\equiv 0,1\pmod{4}$. See Chapter 11 in \cite{A} for more details.

A Wh$(4n+1)$ is {\em $\mathbb{Z}$-cyclic} if the players are taken in $\mathbb{Z}_{4n+1}$ and the round $j+1 \pmod{4n+1}$ is obtained from round $j$ by adding $1\pmod {4n+1}$ to each element. A Wh$(4n)$ is {\em $\mathbb{Z}$-cyclic} if the players are the elements of $\mathbb{Z}_{4n-1}\cup\{\infty\}$ and the rounds are similarly cyclically generated, where $\infty+1=\infty$. Conventionally, in a $\mathbb{Z}$-cyclic Wh$(4n+1)$, 0 is missing from the initial round.

Much less is known about the existence of $\mathbb{Z}$-cyclic whist tournaments despite of the efforts of many authors, for example, Anderson, Finizio and Leonard \cite{AFL}, Buratti \cite{Bur}, Feng and Chang \cite{FC}, Ge and Ling \cite{GLi}, etc. The interested reader is referred to \cite{AFGG} and the references therein.

For any odd positive integer $u$, the set $\{\{x,-x\}:x\in \mathbb{Z}_u\setminus \{0\}\}$ is called the {\em patterned starter} for $\mathbb{Z}_u$, and it is sometimes convenient to call the set $\{\{x,-x\}:x\in \mathbb{Z}_u\setminus \{0\}\}\cup\{\{\infty,0\}\}$ the {\em patterned starter} for $\mathbb{Z}_u\cup \{\infty\}$. A {\em $\mathbb{Z}$-cyclic patterned starter Wh$(v)$}, briefly ZCPS-Wh$(v)$, is a $\mathbb{Z}$-cyclic Wh$(v)$ if the collection of its initial round partner pairs forms the patterned starter for $\mathbb{Z}_v$ when $v=4n+1$ or $\mathbb{Z}_{v-1}\cup \{\infty\}$ when $v=4n$. Finizio \cite{f1994} first introduced the concept of $\mathbb{Z}$-cyclic patterned starter whist tournaments.

\begin{Example}\label{eg:ZCPS-Wh} {\rm \cite[Table 2]{f1994}}
The initial round of a ZCPS-Wh$(13)$ over $\mathbb{Z}_{13}$ is given by the following games: $(1,5,-1,-5)$, $(2,3,-2,-3)$, $(4,6,-4,-6)$.
\end{Example}

\begin{Example}\label{eg:ZCPS-Wh-1} {\rm \cite[Example 1.1]{fl}}
The initial round of a ZCPS-Wh$(28)$ over $\mathbb{Z}_{27}\cup \{\infty\}$ is given by the following games: $(\infty,3,0,-3)$, $(4,5,-4,-5)$, $(6,10,-6,-10)$, $(11,9,-11,-9)$, $(12,7,-12,-7)$, $(8,2,-8,-2)$, $(13,1,-13,-1)$.
\end{Example}


\begin{Proposition}\label{PSandZCPS}
\begin{enumerate}
\item[$(1)$] For $v\equiv 1 \pmod{4}$, there exists a ZCPS-Wh$(v)$ if and only if there exists a PS$(v)$;
\item[$(2)$] For $v\equiv 3 \pmod{4}$, there exists a ZCPS-Wh$(v+1)$ if and only if there exists an APS$(v,\alpha,\alpha)$ for some $\alpha\in \mathbb Z_v$.
\end{enumerate}
\end{Proposition}

\proof Suppose there exists a ZCPS-Wh$(v)$ for $v\equiv 1 \pmod{4}$ (resp. ZCPS-Wh$(v+1)$ for $v\equiv 3 \pmod{4}$) and denote by ${\cal R}_0$ its initial round. Then the games of ${\cal R}_0$ are of the form $(x,y,-x,-y)$, and when $v\equiv 3 \pmod{4}$ we have a special game that has partner pairs $\{0,\infty\}$ and $\{\alpha,-\alpha\}$. Denote by $\cal S$ the set of pairs $\{x,y\}$ of players at the first two positions in the games of ${\cal R}_0$ that do not contain $\infty$. By Condition (2) in the definition of whist tournaments, we have
\begin{eqnarray*}
\begin{array}{ll}
\bigcup\limits_{\{x,y\}\in {\cal S}} \pm\{x,y\}=\mathbb{Z}_v\setminus \{0\} &\mbox{ if } v\equiv 1 \ ({\rm mod}\ 4);\\
\bigcup\limits_{\{x,y\}\in {\cal S}} \pm\{x,y\}=\mathbb{Z}_v\setminus\{0,\alpha,-\alpha\} &\mbox{ if } v\equiv 3\ ({\rm mod}\ 4).
\end{array}
\end{eqnarray*}
By Condition (4) in the definition of whist tournaments, since the tournament is $\mathbb{Z}$-cyclic we have
\begin{eqnarray*}
\begin{array}{ll}
\bigcup\limits_{(x,y)\in {\cal S}} \pm\{x+y,x-y\}=\mathbb{Z}_v\setminus\{0\} &\mbox{ if } v\equiv 1 \ ({\rm mod}\ 4);\\
\bigcup\limits_{(x,y)\in {\cal S}} \pm\{x+y,x-y\}=\mathbb{Z}_v\setminus\{0,\alpha,-\alpha\} &\mbox{ if } v\equiv 3\ ({\rm mod}\ 4).
\end{array}
\end{eqnarray*}
Then $\cal S$ forms a PS$(v)$ when $v\equiv 1 \pmod{4}$ (resp. an APS$(v,\alpha,\alpha)$ when $v\equiv 3 \pmod{4}$). Conversely, reverse the above process to complete the proof. \endproof

\subsubsection{Known results on partitionable sets}

Zhang and Chang \cite{ZC2009} in 2009 showed that any PS$(v)$ with $v\equiv 9\pmod {12}$ and $v\not\equiv 81\pmod{108}$ cannot exist. They also conjectured that no PS$(v)$ with $v\equiv 81\pmod{108}$ exist. Abel, Anderson and Finizio \cite{AAF2011} in 2011 made the same conclusion and the same conjecture for ZCPS-Wh$(v)$. Hu and Ge \cite{HG2012} in 2012 confirmed the conjecture.

\begin{Proposition}\label{prop:nece PS} {\rm (Necessary Condition)}
{\rm \cite[Theorem 2.1]{HG2012}}
Let $v\equiv 1 \pmod{4}$. If there exists a ZCPS-Wh$(v)$, i.e. a PS$(v)$, then $v\equiv 1,5 \pmod{12}$.
\end{Proposition}

Bose and Cameron \cite{BC} in 1965 and Baker \cite{B} in 1975 independently established the existence of ZCPS-Wh$(p)$ for any prime $p\equiv 1 \pmod{4}$. Zhang and Chang \cite{ZC2009} in 2009 presented the same result for PS$(p)$. In fact Watson \cite{W} in 1954 gave a more general result for the existence of ZCPS-Whs.

\begin{Theorem}{\rm \cite{W}}\label{prop:W}
There exists a ZCPS-Wh$(v)$, i.e., a PS$(v)$, if $v$ is a finite product of primes each congruent to $1$ modulo $4$.
\end{Theorem}

Theorem \ref{prop:W} can be obtained by the following standard recursive construction together with the existence of a ZCPS-Wh$(p)$ for any prime $p\equiv 1 \pmod{4}$.

\begin{Construction}{\rm \cite[Theorem 3.2]{AFL}}\label{thm:recur}
Let $u,v\equiv 1\pmod{4}$. Then the existence of a ZCPS-Wh$(u)$ and a ZCPS-Wh$(v)$, i.e., a PS$(u)$ and a PS$(v)$, implies the existence of a ZCPS-Wh$(uv)$, i.e., a PS$(uv)$.
\end{Construction}

We summarize all the other known results on the existence of ZCPS-Wh$(v)$ with $v\equiv 1\pmod{4}$ as follows.

\begin{Theorem}\label{prop:PS known}
There exists a ZCPS-Wh$(v)$, i.e., a PS$(v)$, for each of the following cases:
\begin{enumerate}
\item[$(1)$] {\rm \cite[Lemma 5.6]{HG2013}} $v\leq 300$ and $v\equiv 1,5 \pmod{12}$;
\item[$(2)$] {\rm \cite{LJ2008}} $v=p^2$ where $3< p< 3500$ and $p\equiv 3 \pmod{4}$ is a prime.
\end{enumerate}
\end{Theorem}

\subsubsection{Known results on almost partitionable sets with $\alpha=\beta$}

For ZCPS-Wh$(v+1)$ with $v\equiv 3\pmod{4}$, the known results are rare until very recently. Abel, Anderson and Finizio \cite{AAF2011} developed necessary conditions for the existence of such tournaments. Hu and Ge \cite{HG2012,HG2013} improved their results.

\begin{Proposition}\label{prop:nece APS} {\rm (Necessary Condition)}
\begin{enumerate}
\item[$(1)$] {\rm \cite[Theorem 2.2]{HG2012}} Let $v\equiv 7,11\pmod{12}$. If there exists a ZCPS-Wh$(v+1)$, i.e., an APS$(v,\alpha,\alpha)$ for some $\alpha\in \mathbb Z_v$, then $v$ is not square free.
\item[$(2)$] {\rm \cite[Theorem 2.1]{HG2013}} Let $v\equiv 3\pmod{12}$. If there exists a ZCPS-Wh$(v+1)$, i.e., an APS$(v,\alpha,\alpha)$ for some $\alpha\in \mathbb Z_v$, then $v=3^{2a+1} v_1$, where $a\geq 0$ and $v_1\equiv 1\pmod{12}$.
\end{enumerate}
\end{Proposition}

\begin{Theorem}\label{SmallAPS}{\rm \cite[Lemma 5.11]{HG2013}}
For $v\equiv 3\pmod{4}$ and $v\leq 299$, a ZCPS-Wh$(v+1)$, i.e., an APS$(v,\alpha,\alpha)$ for some $\alpha\in \mathbb Z_v$, exists for any $v$ that satisfies the necessary conditions in Proposition $\ref{prop:nece APS}$ except possibly for $v\in\{243,255,275\}$, i.e., $v\in \{3,27,39,75,111,147,175,183,219,291\}$.
\end{Theorem}

\begin{Theorem}\label{infiAPS}
There exists a ZCPS-Wh$(v+1)$, i.e. an APS$(v,\alpha,\alpha)$ for some $\alpha\in \mathbb Z_v$, for each of the following cases:
\begin{enumerate}
\item[$(1)$] {\rm \cite[Proposition 5.14]{HG2013}} for $v=3p^2$ where $p\equiv 11\pmod{12}$ is a prime and $11\leq p\leq 359$;
\item[$(2)$] {\rm \cite[Theorem 4.9]{HG2012}} $v=3p$ where $p$ is a prime and either $p\equiv 13\pmod{24}$, or $p\equiv 1\pmod{24}$ and $p>1.9\times 10^{12}$ $($or $p\leq 1201)$;
\item[$(3)$] {\rm \cite[Theorem 4.6]{HG2013}} $v=27 p$ where $p\equiv 13\pmod{24}$ is a prime and $p\geq 37$.
\end{enumerate}
\end{Theorem}

\subsubsection{Other variants of whist tournaments}

Ever since the existence of whist tournaments was completely settled, the focus has turned to whist tournaments with additional properties. A game ($a,b,c,d$) can be seen as a cyclic order of the four players sitting round a table.

A {\em directed whist tournament} (DWh) is a whist tournament having the property that every player has every other player as an opponent on his left exactly once
and as an opponent on his right exactly once. A DWh$(v)$ is equivalent to a {\it resolvable $(v, 4, 1)$-perfect Mendelsohn design} \cite{BZ}. By collecting all known results on resolvable perfect Mendelsohn designs, Abel, Bennett and Ge \cite{ABG} showed that a DWh$(v)$ exists if and only if $v \equiv 0,1 \pmod{4}$ except for $v= 4,$ $8,$ $12$ and possibly for $v\in\{16,$ $20,$ $24,$ $32,$ $36,$ $44,$
$48,$ $52,$ $56,$ $64,$ $68,$ $76,$ $84,$ $88,$ $92,$ $96,$ $104,$ $108,$
$116,$ $124,$ $132,$ $148,$ $152,$ $156,$ $172,$ $184,$ $188\}$.

Baker and Wilson \cite{BW} pointed out that no $\mathbb{Z}$-cyclic DWh$(v)$ exists whenever $v\equiv 0\pmod{4}$, but no proof was provided. Finizio and Leonard gave a proof in \cite{fl}.

\begin{Proposition}\label{no DWh}{\rm \cite[Theorem 3.8]{fl}}
If $v\equiv 0\pmod{4}$, then there is no $\mathbb{Z}$-cyclic DWh$(v)$.
\end{Proposition}

If there exists a $\mathbb{Z}$-cyclic DWh$(4n+1)$ with games $(a_i,b_i,c_i,d_i)$, $1\leq i\leq n$, in the initial round, then by the ``directed'' property, we have
$$\bigcup_{i=1}^n \{b_i-a_i,c_i-b_i,d_i-c_i,a_i-d_i\}=\mathbb Z_v\setminus \{0\}.$$
Therefore, the following proposition holds.

\begin{Proposition}\label{ZCPS to DWh}{\rm \cite[Theorem 3.5]{fl}}
Every ZCPS-Wh$(4n+1)$, i.e., a PS$(4n+1)$, is a $\mathbb{Z}$-cyclic DWh$(4n+1)$.
\end{Proposition}

In the game ($a,b,c,d$), $a$ and $c$ are said to be {\em partners of the first kind}, while $b$ and $d$ are said to be {\em partners of the second kind}. An {\em ordered whist tournament} (OWh) is a whist tournament in which each player opposes every other player exactly once as a partner of the first kind and exactly once as a partner of the second kind. It is known that an OWh$(v)$ exists if and only if $v\equiv 1\pmod{4}$ \cite{ACF,cf}.

If there exists a $\mathbb{Z}$-cyclic OWh$(4n+1)$ with games $(a_i,b_i,c_i,d_i)$, $1\leq i\leq n$, in the initial round, then by the ``ordered'' property, we have
$$\bigcup_{i=1}^n \{b_i-a_i,d_i-a_i,b_i-c_i,d_i-c_i\}=\mathbb Z_v\setminus \{0\}.$$
Therefore, the following proposition holds.

\begin{Proposition}\label{ZCPS to OWh}{\rm \cite[Theorem 4.5]{GZ}}
Every ZCPS-Wh$(4n+1)$, i.e., a PS$(4n+1)$, is a $\mathbb{Z}$-cyclic OWh$(4n+1)$.
\end{Proposition}

We remark that Abel, Costa and Finizio \cite{ACF} also stated Propositions \ref{ZCPS to DWh} and \ref{ZCPS to OWh} as their Corollary 5.

\subsection{Difference matrices}\label{sec:DM}

A {\em cyclic $(v,k,1)$-difference matrix}, briefly CDM, is a $k\times v$ array $D=(d_{ij})$ over $\mathbb Z_v$ satisfying that for each $1\leq r<s\leq k$, $\{ d_{rl}-d_{sl}: 1\leq l\leq v \}=\mathbb Z_v$.

If there exists a $\mathbb{Z}$-cyclic DWh$(4n+1)$ with games $(a_i,b_i,c_i,d_i)$, $1\leq i\leq n$, in the initial round, then let
\begin{center}
$A_i=\left(
\begin{array}{cccc}
0 & 0 & 0 & 0 \\
a_i & b_i & c_i & d_i \\
b_i & c_i & d_i & a_i \\
c_i & d_i & a_i & b_i \\
d_i & a_i & b_i & c_i \\
\end{array}
\right)
$.
\end{center}
It is readily checked that $[O|A_1|A_2|\cdots | A_n]$ forms a CDM$(4n+1,5,1)$, where $O=(0,0,0,0,0)^T$ is a column vector.

\begin{Proposition}\label{DWh to DM}{\rm \cite[Corollary 6]{f}}
If there exists a $\mathbb{Z}$-cyclic DWh$(4n+1)$, then there exists a CDM$(4n+1,5,1)$.
\end{Proposition}

Combining Propositions \ref{ZCPS to DWh} and \ref{DWh to DM}, we have the following corollary.

\begin{Corollary}\label{cor:CDM}
If there exists a ZCPS-Wh$(4n+1)$, i.e., a PS$(4n+1)$, then there exists a CDM$(4n+1,5,1)$.
\end{Corollary}

\subsection{Cyclic balanced sampling plans excluding contiguous units}\label{sec:BSEC}

The existence problem of partitionable sets was initially studied in \cite{ZC2009} and, as an application, several infinite families of cyclic BSECs with block size four were presented.

Let $X=(x_0,x_1,\ldots,x_{v-1})$ be cyclically ordered. For $0\leq i \leq v-1$, $x_i$ and $x_{i+1}$ are said to be {\em contiguous points}, where the addition is reduced modulo $v$.
A BSEC$(v,k,1)$ is a pair $(X,\mathcal{B})$ where $X$ is a set of $v$ points in a cyclic ordering and $\mathcal{B}$ is a collection of $k$-subsets of $X$ called {\em blocks}, such that any two contiguous points do not appear in any block while any two non-contiguous points appear in exactly one block.
If a BSEC$(v,k,1)$ admits $\mathbb{Z}_v$ as its automorphism group, then it is called {\em cyclic} and written simply as a CBSEC$(v,k,1)$.

Zhang and Chang \cite{ZC2009} showed that if there exists a PS$(v)$ for $v\equiv 1,5 \pmod{12}$, then there exists a CBSEC$(3v,4,1)$. By Proposition \ref{prop:nece PS}, a PS$(v)$ exists only if $v\equiv 1,5 \pmod{12}$. Therefore, the following theorem is stated.

\begin{Theorem}{\rm \cite[Constructions 3.2]{ZC2009}}
If there is a PS$(v)$, then there is a CBSEC$(3v,4,1)$.
\end{Theorem}

Construction 3.3 in \cite{ZC2009} provides a recursive construction for CBSECs: if there exist a PS$(v)$, a CBSEC$(u,4,1)$ and a $(v,4,1)$-CDM, then there exists a CBSEC$(uv,4,1)$. By Corollary \ref{cor:CDM}, this construction can be simplified as follows.

\begin{Theorem}
If there exist a PS$(v)$ and a CBSEC$(u,4,1)$, then there exists a CBSEC$(uv,4,1)$.
\end{Theorem}

\subsection{Disjoint difference families}

Let $(G,+)$ be a group of order $v$. A $(G,k,\lambda)$ {\em difference family} is a family $\mathfrak{B}=[B_1,\dots,B_s]$ of $k$-subsets of $G$ such that the list
$\bigcup_{i=1}^s[x-y:x,y\in B_i, x\not=y]$ covers every element of $G\setminus \{0\}$ exactly $\lambda$ times. Clearly $s=\lambda(v-1)/(k(k-1))$. A $(G,k,\lambda)$-difference family is said to be {\em disjoint} and written as a $(G,k,\lambda)$-DDF when its base blocks are mutually disjoint.

Similar arguments to those in the proof of Proposition \ref{PSandZCPS}(1) show that the following equivalence holds between partitionable sets and a special kind of disjoint difference families.

\begin{Proposition}\label{prop:PS-DDF}
A PS$(G)$ in an abelian group $(G,+)$ is equivalent to a $(G,4,3)$-DDF whose base blocks are all of the form $\{x,-x,y,-y\}$.
\end{Proposition}

By exploring Ferrero pairs, Buratti \cite{b19} presented a simple proof for the existence of $(G,k,k-1)$-DDFs where $G$ is the additive group of $\mathbb{F}_{q_1}\times\cdots\times\mathbb{F}_{q_r}$ with each $q_i\equiv 1\pmod{k}$ a prime power, and established the existence of $(G,k,k-1)$-DDFs for any abelian group $G$ of order $v$ where  every prime factor of $v$ is congruent to $1 \pmod{k}$. Taking $k=4$ and examining the proof of these two infinite families of DDFs in \cite[Corollaries 3.3 and 3.5]{b19}, one can see that each base block of the resulting DDFs is of the form $\{x,-x,y,-y\}$. Thus the following theorem can be obtained.

\begin{Theorem}\label{thm:abe-1}
\begin{enumerate}
\item[$(1)$]
For $1\leq i\leq r$, let $q_i\equiv 1\pmod{4}$ be a prime power. Then there exists a PS$(G)$ where $G$ is the additive group of $\mathbb{F}_{q_1}\times\cdots\times\mathbb{F}_{q_r}$.
\item[$(2)$]
If every prime factor of $v$ is congruent to $1 \pmod{4}$, then there exists a PS$(G)$ for any abelian group $G$ of order $v$.
\end{enumerate}
\end{Theorem}


\section{Necessary conditions for almost partitionable sets}

Proposition \ref{prop:nece APS} shows necessary conditions for the existence of a ZCPS-Wh$(v+1)$, i.e., an APS$(v,\alpha,\alpha)$ for some $\alpha\in \mathbb Z_v$. In this section we provide necessary conditions for the existence of an APS$(v,\alpha,\beta)$. 

\begin{Lemma}\label{nec}
If there exists an APS$(v,\alpha,\beta)$, then
\begin{eqnarray*}
2\alpha^2-\beta^2\equiv
\left\{
\begin{array}{ll}
    \frac{v}{3} \ ({\rm mod }\ v), & \hbox{ \rm{if} $v\equiv 3\ ({\rm mod }\ 12)$,}\\[0.4em]
	0\ ({\rm mod }\ v), & \hbox{ \rm{if} $v\equiv 7,11 \ ({\rm mod }\ 12)$.}
\end{array}
\right.
\end{eqnarray*}
\end{Lemma}

\proof Denote by $\cal S$ an APS$(v,\alpha,\beta)$. By the definition of APS, $\bigcup\limits_{\{x,y\}\in {\cal S}} \pm\{x,y\}=\mathbb{Z}_v\setminus \{0,\alpha,-\alpha\}$ and $\bigcup\limits_{\{x,y\}\in {\cal S}} \pm\{x-y,x+y\}=\mathbb{Z}_v\setminus \{0,\beta,-\beta\}$. It follows that we can compute the sum of squares of all elements in $\mathbb{Z}_v$ in three ways as follows:
\begin{eqnarray*}
2\alpha^2+2\sum_{\{x,y\}\in {\cal S}}(x^2+y^2)\equiv 2\beta^2+4\sum_{\{x,y\}\in {\cal S}}(x^2+y^2)\equiv \sum_{j=1}^{v-1}j^2\ ({\rm mod}\ v).
\end{eqnarray*}
Since $\sum\limits_{j=1}^{v-1}j^2 \equiv 2\sum\limits_{j=1}^{v-1}j^2-\sum\limits_{j=1}^{v-1}j^2\ ({\rm mod}\ v)$, replacing the first $\sum\limits_{j=1}^{v-1}j^2$ in the right hand side by $2\alpha^2+2\sum_{\{x,y\}\in {\cal S}}(x^2+y^2)$ and the second $\sum\limits_{j=1}^{v-1}j^2$ by $2\beta^2+4\sum_{\{x,y\}\in {\cal S}}(x^2+y^2)$, we have
\begin{align*}
\sum\limits_{j=1}^{v-1}j^2 \equiv  2\left(2\alpha^2+2\sum\limits_{\{x,y\}\in S}(x^2+y^2)\right)-\left(2\beta^2+4\sum\limits_{\{x,y\}\in S}(x^2+y^2)\right) \equiv 4\alpha^2-2\beta^2\ ({\rm mod}\ v).
\end{align*}
On the other hand,
\begin{align*}
\sum\limits_{j=1}^{v-1}j^2 = \frac{(v-1)v(2v-1)}{6} \equiv\left\{
\begin{array}{lll}
\frac{2v}{3}\ ({\rm mod}\ v), & {\rm if}\ v\equiv 3\ ({\rm mod}\ 12); \\
0\ ({\rm mod}\ v), & {\rm if}\ v\equiv 7,11\ ({\rm mod}\ 12).
\end{array}
\right.
\end{align*}
Thus
\begin{align*}
4\alpha^2-2\beta^2 \equiv\left\{
\begin{array}{lll}
\frac{2v}{3}\ ({\rm mod}\ v), & {\rm if}\ v\equiv 3\ ({\rm mod}\ 12); \\
0\ ({\rm mod}\ v), & {\rm if}\ v\equiv 7,11\ ({\rm mod}\ 12).
\end{array}
\right.
\end{align*}
Due to $\gcd(2,v)=1$, the conclusion follows immediately. \endproof

\begin{Corollary}\label{cor:nec}
Let $v\equiv 3\ ({\rm mod}\ 12)$ and $v=3p_1p_2\cdots p_s$, where the $p_i$'s are different primes such that $p_i\equiv \pm3\ ({\rm mod}\ 8)$ for each $i\in\{1,2,\ldots,s\}$. If there exists an APS$(v,\alpha,\beta)$, then $\alpha,\beta\in\{v/3,2v/3\}$.
\end{Corollary}

\proof By Lemma \ref{nec}, $v\equiv 3\ ({\rm mod}\ 12)$ yields $2\alpha^2\equiv \beta^2 \pmod{\frac{v}{3}}$. Since $2$ is not a square in ${\mathbb Z}^*_{p}$ for any prime $p\equiv \pm3 \pmod{8}$, we have that $\alpha\equiv \beta \equiv 0 \pmod{\frac{v}{3}}$. \endproof

As expected, Lemma \ref{nec} provides weaker necessary conditions for the existence of APS$(v,\alpha,\beta)$ than those for the existence of APS$(v,\alpha,\alpha)$ in Proposition \ref{prop:nece APS}. For example, by Proposition \ref{prop:nece APS}(1), an APS$(7,\alpha,\alpha)$ does not exist for any $\alpha\in \mathbb Z_7$, but we can construct an APS$(7,2,1)$ consisting of only one pair $\{1,4\}$, which satisfies Lemma \ref{nec}.

Now we apply Lemma \ref{nec} to determine for what values of $v$, an APS$(v,\alpha,\beta)$ for any $\alpha$ and $\beta$ cannot exist.

\begin{Lemma}\label{Nec}
In each of the following cases, there is no APS$(v,\alpha,\beta)$ for any $\alpha$ and $\beta$:
\begin{enumerate}
\item[$(1)$] $v\equiv 3\ ({\rm mod}\ 12)$ and $v=3^{2a}v_1$, where $a\geq 1$ and $3\nmid v_1;$
\item[$(2)$] $v\equiv 15\ ({\rm mod}\ 36)$ and $v=3p_1p_2\cdots p_s$, where the $p_i$'s are different primes such that $p_i\equiv \pm3\ ({\rm mod}\ 8)$ for each $i\in\{1,2,\ldots,s\};$
\item[$(3)$] $v\equiv 7,11\ ({\rm mod}\ 12)$ and $v=p_1p_2\cdots p_t$, where the $p_i$'s are different primes such that $p_i\equiv \pm3\ ({\rm mod}\ 8)$ for each $i\in\{1,2,\ldots,t\}$.
\end{enumerate}
\end{Lemma}

\proof $(1)$ By Lemma \ref{nec}, it suffices to show that in the given assumption of $v$, for any $\alpha,\beta\in \mathbb{Z}_v\setminus\{0\}$, $2\alpha^2-\beta^2\neq m\cdot\frac{v}{3}$, where $m\equiv 1\ ({\rm mod}\ 3)$. Clearly the maximal power of $3$ that divides $2\alpha^2-\beta^2$ is even, while the maximal power of $3$ that divides $m\cdot\frac{v}{3}=3^{2a-1}\cdot mv_1$ is odd since $\gcd(mv_1,3)=1$. Hence $2\alpha^2-\beta^2\neq m\cdot\frac{v}{3}$.

$(2)$ $v\equiv 15\ ({\rm mod}\ 36)$ implies that $w:= v/3 \equiv 2\ ({\rm mod}\ 3)$. By Corollary \ref{cor:nec} we have $\alpha=wa$ and $\beta=wb$ for suitable $a,b\in\{1,2\}$. By Lemma \ref{nec}, if there exists an APS$(v,\alpha,\beta)$, then $2\alpha^2-\beta^2\equiv w\ ({\rm mod }\ 3w)$. This gives $w^2(2a^2-b^2)\equiv w\ ({\rm mod}\ 3w)$. Reducing modulo $3$ we get $1 \equiv 2\ ({\rm mod}\ 3)$ that is absurd.

$(3)$ By Lemma \ref{nec}, it suffices to show that in the given assumption of $v$, for any $\alpha,\beta\in \mathbb{Z}_v\setminus \{0\}$, $2\alpha^2-\beta^2\not\equiv 0\ ({\rm mod}\ v)$. This is straightforward since $2$ is not a square in ${\mathbb Z}^*_{p}$ for any prime $p\equiv \pm3 \pmod{8}$. \endproof

Next we show that, for any $v\equiv 3\pmod {4}$ that does not satisfy the conditions $(1)$, $(2)$ and $(3)$ in Lemma \ref{Nec}, there exists a pair of $\alpha$ and $\beta$ satisfying the condition in Lemma \ref{nec}.

\begin{Proposition}
Let $v\equiv 3\pmod {4}$ and $v$ does not satisfy the conditions $(1)$, $(2)$ and $(3)$ in Lemma $\ref{Nec}$.
\begin{itemize}
\item[$(1)$] If $v\equiv 3\ ({\rm mod}\ 12)$, then there exist nonzero $($modulo $v)$ $\alpha$ and $\beta$ such that $2\alpha^2-\beta^2\equiv \frac{v}{3}\pmod{v}.$
\item[$(2)$] If $v\equiv 7,11\ ({\rm mod}\ 12)$, then  there exist nonzero $($modulo $v)$ $\alpha$ and $\beta$ such that $2\alpha^2-\beta^2\equiv 0\pmod{v}$.
\end{itemize}
\end{Proposition}

\proof (1) Let $v\equiv 3 \pmod{12}$ that does not satisfy the conditions of Lemma $\ref{Nec}$. Then $v=3^dp_1^{\gamma_1}p_2^{\gamma_2}\cdots p_s^{\gamma_s}$ for some odd integer $d$, where the $p_i$'s are different primes not equal to $3$. There are the following possibilities:
\begin{center}
\begin{itemize}
\item[A)] $p_1^{\gamma_1}p_2^{\gamma_2}\cdots p_s^{\gamma_s}\equiv 1\ ({\rm mod}\ 12)$, which implies $v/3\equiv 3^{d-1}\pmod{3^d}$;
\item[B)] $p_1^{\gamma_1}p_2^{\gamma_2}\cdots p_s^{\gamma_s}\equiv 5\ ({\rm mod}\ 12)$, which implies $v/3\equiv -3^{d-1}\pmod{3^d}$,
\begin{itemize}
\item[B.1)] $d>1$;
\item[B.2)] $d=1$.
\end{itemize}
\end{itemize}
\end{center}
By the Chinese remainder theorem, it suffices to find $\alpha_0,\beta_0$ such that $2\alpha_0^2-\beta_0^2\equiv v/3\pmod{3^d}$ and $\alpha_i,\beta_i$ such that $2\alpha_i^2-\beta_i^2\equiv 0\pmod{p_i^{\gamma_i}}$ for $1\leq i\leq s$, where at least one $\alpha_j$ is nonzero and one $\beta_{j'}$ is nonzero.

For Cases A and B.1, take $\alpha_i\equiv\beta_i\equiv 0\pmod{p_i^{\gamma_i}}$ for $1\leq i\leq s$. In Case A, since $v/3\equiv 3^{d-1}\pmod{3^d}$, we take $\alpha_0\equiv\beta_0\equiv\sqrt{3^{d-1}}\pmod{3^d}$ (note that $d$ is odd); in Case B.1, since $d>1$ and $v/3\equiv -3^{d-1}\pmod{3^d}$, we take $\alpha_0\equiv 3^{d-1}\pmod{3^d}$ and $\beta_0\equiv \sqrt{3^{d-1}}\pmod{3^d}$.

For Case B.2, by Lemma \ref{Nec}(2), without loss of generality we can assume that either $\gamma_1>1$ or $p_1\equiv \pm1 \pmod{8}$. Take $\alpha_0\equiv 0\pmod{3}$ and $\beta_0\equiv 1\pmod{3}$ such that $2\alpha_0^2-\beta_0^2\equiv v/3\equiv-1 \pmod{3}$. Take $\alpha_i\equiv\beta_i\equiv 0\pmod{p_i^{\gamma_i}}$ for $2\leq i\leq s$ when $s\geq 2$. If $\gamma_1>1$, we take $\alpha_1\equiv \beta_1\equiv p_1^{\gamma_1-1}\pmod{p_1^{\gamma_1}}$; if $p_1\equiv \pm1 \pmod{8}$, since $2$ is a square in $\mathbb Z^*_{p_1}$, there exist $\alpha_1, \beta_1\not\equiv 0\pmod{p_1^{\gamma_1}}$ such that $2\alpha_1^2-\beta_1^2\equiv 0 \pmod{p_1^{\gamma_1}}$.

(2) Let $v\equiv 7,11 \pmod{12}$ that does not satisfy the condition of Lemma $\ref{Nec}$. Write $v=p_1^{\gamma_1}p_2^{\gamma_2}\cdots p_t^{\gamma_t}$, where the $p_i$'s are different primes. Without loss of generality we can assume that either $\gamma_1>1$ or $p_1\equiv \pm1 \pmod{8}$. Then the same argument as that in Case B.2 shows that there exist $\alpha_1, \beta_1\not\equiv 0\pmod{p_1^{\gamma_1}}$ such that $2\alpha_1^2-\beta_1^2\equiv 0 \pmod{p_1^{\gamma_1}}$. Then the proof is completed by using the Chinese remainder theorem.
\endproof

\section{A silver ratio construction for almost partitionable sets}

Let $\mathbb{Z}_p^*$ be the multiplicative group of nonzero integers modulo a prime $p$. For any prime $p\equiv 1,7 \pmod{8}$, $2$ is a square in $\mathbb{Z}_p^*$. The {\em silver elements} of $\mathbb{Z}_p^*$ are defined to be the two elements $1+x$ and $1-x$ of $\mathbb{Z}_p^*$ where $x$ and $-x$ are the square roots of $2$ modulo $p$. Usually we write $\sqrt{2}$ to represent one of these roots.
Buratti introduced silver elements to construct optical orthogonal codes in \cite{b18}. The proof of Theorem 3.1 in \cite{b18} uses an APS$(p,1,\sqrt{2})$ implicitly, and so the following lemma has been essentially proved in \cite{b18}.

\begin{Lemma}\label{Marco}
Let $p\equiv 7 \pmod{8}$ be a prime. If $\theta=1+\sqrt{2}$ generates $\mathbb{Z}_p^*/\{1,-1\}$, then there exists an APS$(p,1,\theta-1)$.
\end{Lemma}

\proof Consider the set $\cal S$ of unordered pairs from  $\mathbb{Z}_{p}$ given by
$${\cal S}=\left\{\{\theta^{2i-1},\theta^{2i}\}\mid 1\leq i\leq \frac{p-3}{4}\right\}.$$
Since $\theta$ satisfies the equation $x(x-1)=x+1$, we have
$$\bigcup_{i=1}^{(p-3)/4}\pm\left\{\theta^{2i-1},\theta^{2i}\right\}=\mathbb{Z}_{p}\setminus \{0,1,-1\},$$
and
$$\bigcup_{i=1}^{(p-3)/4}\pm\{\theta^{2i-1}(\theta-1),\theta^{2i-1}(\theta+1)\}=
\bigcup_{i=1}^{(p-3)/4}\pm\{\theta^{2i-1}(\theta-1),\theta^{2i}(\theta-1)\}=\mathbb{Z}_{p}\setminus \{0,\theta-1,1-\theta\}.$$
Therefore $\cal S$ is an APS$(p,1,\theta-1)$.
\endproof

\begin{Theorem}\label{cor:Marco}
Suppose that $p\equiv 7 \pmod{8}$ is a prime, and $1+\sqrt{2}$ generates $\mathbb{Z}_p^*/\{1,-1\}$. Then an APS$(p,\alpha,\beta)$ exists if and only if $2\alpha^2-\beta^2\equiv0\pmod{p}$.
\end{Theorem}

\proof By Lemma \ref{nec}, an APS$(p,\alpha,\beta)$ with $p\equiv 7 \pmod{8}$ a prime exists only if $\beta\equiv \pm \sqrt{2}\alpha\pmod{p}$. Let $\cal S$ be an APS$(p,1,\sqrt{2})$ given by Lemma \ref{Marco}. Then
$\{\{\alpha x,\alpha y\} \mid\{x,y\}\in{\cal S}\}$
is an APS$(p,\alpha,\beta)$. \endproof

One can check that when $p\in \{7,23,31,47,71,127,151,167,191,263,271\}$, $1+\sqrt{2}$ generates $\mathbb{Z}_p^*/\{1,-1\}$. Thus we have the following corollary.

\begin{Corollary}\label{cor:APS from Marco}
There is an APS$(p,\alpha,\beta)$ for any $p\in \{7,23,31,47,71,127,151,167,191,263,271\}$ and $2\alpha^2-\beta^2\equiv0\pmod{p}$.
\end{Corollary}

The APS given by Lemma \ref{Marco} has all its pairs that are multiples of a ``starter pair". Are there other APSs with this property? We shall show that when $p\equiv 3 \pmod{4}$ is a prime, any APS$(p)$ whose pairs are all multiples of a ``starter pair" is necessarily the APS$(p,1,\sqrt{2})$ given in Lemma \ref{Marco} possibly multiplied by a nonzero element in $\mathbb{Z}_p$.


\begin{Proposition}\label{prop:new}
Let $p\equiv 3 \pmod{4}$ be a prime and $x\in\mathbb{Z}_p$. Then there exists a set $S\subseteq \mathbb{Z}_p$ such that $\{\{s, sx\}\mid s\in S\}$ is an APS$(p)$ if and only if $x$ is a generator of $\mathbb{Z}^*_p/\{1,-1\}$ and $x^2\pm 2x-1=0$. Furthermore, this APS must be an APS$(p,1,\sqrt{2})$ multiplied by a nonzero element in $\mathbb{Z}_p$.
\end{Proposition}

\proof To prove the sufficiency, let ${\cal A}=\{\{s, sx\}\mid s\in S\}$, where $S=\left\{x^{2i-1}\mid 1\leq i\leq (p-3)/4\right\}.$
Similar arguments to those in the proof of Lemma \ref{Marco} show that $\cal A$ forms an APS$(p,1,x-1)$ or APS$(p,1,x+1)$ according to $x^2-2x-1=0$ or $x^2+2x-1=0$, respectively. Note that $x^2-2x-1=0$ implies $x-1=\pm \sqrt{2}$ and $x^2+2x-1=0$ implies $x+1=\pm \sqrt{2}$. Thus $\cal A$ is an APS$(p,1,\sqrt{2})$.

Next we examine the necessity. Assume that $\{\{s,sx\}\mid s\in S\}$ is an APS$(p)$ with $p=4n+3$ a prime. Let $g$ be a generator of $\mathbb{Z}^*_p$ and consider the $\log$ map $g^i\in\mathbb{Z}^*_p\rightarrow i\in\mathbb{Z}_{4n+2}$. Set $i=\log(x)$, $j=\log(\frac{x+1}{x-1})$ and $K=\{\log(s)\mid s\in S\}$. Then $\{\{0,i,2n+1,2n+1+i\}+k\mid k\in K\}$ and $\{\{0,j,2n+1,2n+1+j\}+k\mid k\in K\}$ are both sets of disjoint quadruples in $\mathbb{Z}_{4n+2}$. Hence the sets of pairs ${\cal M}_i=\{\{0,i\}+k\mid k\in K\}$ and ${\cal M}_j=\{\{0,j\}+k\mid k\in K\}$ are near-perfect matchings of the complete graph on $\mathbb{Z}_{2n+1}$. By Corollary 2.7 in \cite{b96}, $i$ is a generator of $\mathbb{Z}_{2n+1}$; thus every element of $\mathbb{Z}_{2n+1}$ can be seen as a multiple of $i$. If $mi$ is the missing vertex of ${\cal M}_i$, we clearly have $K=\{mi+(2h+1)i \mid 0\leq h\leq n-1\}$. Given that $mi+i\in K$, we have $mi+i+j\not\in K$ that means $j\not\in \{2hi\mid 0\leq h\leq n-1\}$. Analogously, given that $mi+(2n-1)i=mi-2i\in K$, we have $mi-2i+j\not\in K$, i.e., $j\not\in \{(2h+3)i\mid 0\leq h\leq n-1\}$. We conclude that $j$ is either $i$ or $-i$. This, turning back to $\mathbb{Z}^*_p$, implies that $x$ is a generator of $\mathbb{Z}^*_p/\{1,-1\}$ satisfying the identity $x^{\pm 1}=\pm\frac{x+1}{x-1}$ so that either $x^2=-1$ - that is a contradiction - or $x^2\pm 2x-1=0$. \endproof

\section{Recursive constructions for almost partitionable sets}\label{RecSection}

In this section we shall present several recursive constructions for APSs. First we introduce a concept that can be seen as a generalization of PSs and APSs. A PS (resp. APS) for a group $(G,+)$ is a partition of $G$ except for the zero element (resp. three elements including 0) into ordered pairs such that some properties are satisfied. Naturally, we can define a ``partial'' partitionable set to be a partition of $G$ except for some given elements into ordered pairs such that additional properties are satisfied. This concept will be used to construct optical orthogonal codes in Section 7.2.

More formally, a {\em partial partitionable set} in a group $(G,+)$ with two subsets $A_1$ and $A_2$ of $G$ such that $|A_1|=|A_2|$, denoted by a PPS$(G,A_1,A_2)$, is a set $\cal S$ of $(|G|-|A_1|)/4$ ordered pairs from $G$ such that
\begin{enumerate}
\item[$(1)$] $\bigcup\limits_{(x,y)\in {\cal S}} \pm\{x,y\}=G\setminus A_1;$
\item[$(2)$] $\bigcup\limits_{(x,y)\in {\cal S}} \pm\{x-y,x+y\}=G\setminus A_2$.
\end{enumerate}
When $G$ is an abelian group, we write $\{x,y\}$ instead of $(x,y)$. A PPS$(G,\{0\},\{0\})$ is just a PS$(G)$, and a PPS$(G,\{0,\alpha,-\alpha\},\{0,\beta,-\beta\})$ is just an APS$(G,\alpha,\beta)$.

According to the definition of partial partitionable sets, the following construction is straightforward.

\begin{Construction}\label{con:filling} $($Filling Construction$)$
Let $H$ be a subgroup of $G$ and $A_1$, $A_2\subseteq H$ with $|A_1|=|A_2|$. If there exist a PPS$(G,H,H)$ and a PPS$(H,A_1,A_2)$, then there exists a PPS$(G,A_1,A_2)$.
\end{Construction}

If $v=hn$, then $\mathbb Z_v$ has a subgroup $U$ of order $h$. A PPS$(\mathbb Z_v,U,U)$ is often called a {\em $\mathbb{Z}$-cyclic patterned starter whist tournament frame}, denoted by a ZCPS-Wh-frame$(h^n)$. Construction \ref{con:filling} unifies and generalizes Lemmas 3.2, 3.3 and 3.4 in \cite{HG2012}, which are referred to as frame constructions for ZCPS-Whs.

\begin{Construction}\label{con:recursive} $($Product Construction$)$
Let $u$ be an integer such that $\gcd(u,6)=1$. If there is a PPS$(\mathbb Z_{v},A_1,A_2)$, then there is a PPS$(\mathbb{Z}_{vu},B_1,B_2)$, where $B_i=\{a+vs\mid a\in A_i,s\in\mathbb{Z}_{u} \}$ for $i\in\{1,2\}$.
\end{Construction}

\proof Let $\cal T$ be a PPS$(\mathbb Z_{v},A_1,A_2)$. Consider the set $\cal S$ of unordered pairs from  $\mathbb{Z}_{vu}$ given by
$${\cal S}=\left\{\{x+sv,y+2sv\}\mid \{x,y\}\in {\cal T}, s\in \mathbb{Z}_{u}\right\}.$$
Since $\gcd(u,6)=1$, we have
$$\bigcup_{\substack{\{x,y\}\in {\cal T}\\  s\in \mathbb{Z}_{u}}}\pm\{x+sv,y+2sv\}=\mathbb{Z}_{vu}\setminus B_1,$$
and
$$\bigcup_{\substack{\{x,y\}\in {\cal T}\\  s\in \mathbb{Z}_{u}}}\pm\{x+y+3sv,x-y-sv\}=\mathbb{Z}_{vu}\setminus B_2.$$
Therefore $\cal S$ is a PPS$(\mathbb{Z}_{vu},B_1,B_2)$. \endproof

\begin{Corollary}\label{cor:recursive}
If there exist a PS$(v)$ and an APS$(u,\alpha,\beta)$ with $u\equiv 7,11\pmod{12}$, then there exists an APS$(vu,v\alpha,v\beta)$.
\end{Corollary}

\proof Take a PS$(v)$, i.e., a PPS$(\mathbb{Z}_{v},\{0\},\{0\})$. Since $\gcd(u,6)=1$, apply Construction \ref{con:recursive} to obtain a PPS$(\mathbb{Z}_{vu},B,B)$, where $B=\{vs\mid s\in\mathbb{Z}_{u} \}$ is a subgroup of $\mathbb{Z}_{vu}$ and $|B|=u$. Now apply Construction \ref{con:filling} with an APS$(u,\alpha,\beta)$ to get an APS$(vu,v\alpha,v\beta)$. \endproof

The use of Corollary \ref{cor:recursive} relies on the existence of some APS$(u,\alpha,\beta)$s with $u\equiv 7,11\pmod{12}$. We remark that except for the APSs in Theorem \ref{cor:Marco}, by Theorems \ref{SmallAPS} and \ref{infiAPS} the only known APS with such parameters is APS$(175,35,35)$ (see Example 5.7 in \cite{HG2013}).

\begin{Construction}\label{con:recursive-alpha-beta}
Let $p\equiv 7 \pmod{8}$ and $q\equiv 5 \pmod{8}$ be both primes. If there exists an APS$(p,\alpha,\beta)$ for any $\alpha$ and $\beta$ such that $2\alpha^2-\beta^2\equiv0\pmod{p}$, then there exists an APS$(pq,\alpha_1,\beta_1)$ for any $\alpha_1$ and $\beta_1$ such that $2\alpha_1^2-\beta_1^2\equiv0\pmod{pq}$.
\end{Construction}

\proof By Theorem \ref{prop:W}, there exists a PS$(q)$ for any prime $q\equiv 5 \pmod{8}$. A prime $p\equiv 7 \pmod{8}$ implies $p\equiv 7,23 \pmod{24}$. So we can apply Corollary \ref{cor:recursive} with an APS$(p,\alpha,\beta)$ and a PS$(q)$ to obtain an APS$(pq,q\alpha,q\beta)$. If we can show that the number of solutions modulo $p$ of the congruence $2\alpha^2-\beta^2\equiv0\pmod{p}$ is the same as the number of solutions modulo $pq$ of the congruence $2\alpha_1^2-\beta_1^2\equiv0\pmod{pq}$, then the proof is completed.

Since $p\equiv 7 \pmod{8}$ is a prime, $2\alpha^2-\beta^2\equiv0\pmod{p}$ implies $\beta\equiv \pm\sqrt{2}\alpha\pmod{p}$. Thus the number of solutions modulo $p$ of the congruence $2\alpha^2-\beta^2\equiv0\pmod{p}$ is $2p-1$. On the other hand, due to $q\equiv 5 \pmod{8}$, $2$ is not a square modulo $q$. So $2\alpha_1^2-\beta_1^2\equiv0\pmod{q}$ has the unique solution $(\alpha_1,\beta_1)\equiv(0,0)\pmod{q}$. Therefore, by the Chinese remainder theorem, the number of solutions modulo $pq$ of the congruence $2\alpha_1^2-\beta_1^2\equiv0\pmod{pq}$ is also $2p-1$. \endproof

Combining Corollary \ref{cor:APS from Marco} and Construction \ref{con:recursive-alpha-beta}, we have the following result.

\begin{Corollary}\label{cor:APS from recur}
Let $p\in \{7,23,31,47,71,127,151,167,191,263,271\}$ and $q\equiv 5 \pmod{8}$ be a prime. Then there is an APS$(pq,\alpha,\beta)$ for any $\alpha$ and $\beta$ such that $2\alpha^2-\beta^2\equiv0\pmod{pq}$.
\end{Corollary}

\section{Kramer-Mesner method}\label{KMSection}

A celebrated technique for the construction of designs with prescribed automorphism groups is the Kramer-Mesner method (see \cite{KM76} and also \cite{KO06}), which reduces the existence problem of designs to the solution problem of suitable linear systems $MX=J$. Group actions are also used in order to further simplify the problem. The purpose of this section is to provide a similar method to construct  partitionable sets, almost partitionable sets and, more in general, partial partitionable sets.

First of all we introduce some notations about the group action on a set. Given an action of a group $H$ on the set $X$ we will denote by $X/H$ the set of orbits of this action and, for $x\in X$, by $[x]$ the orbit of the element $x$ that is $[x]:=\{x^h\mid h\in H\}$ where $x^h$ is the image of $x$ under $h$. Let ${X}\choose{2}$ be the set of all unordered pairs of $X$. Similarly we will denote by ${{X}\choose{2}}/H$ the set of orbits of the action of $H$ on ${X}\choose{2}$ and, for $B=\{x,y\}\in {{X}\choose{2}}$, by $[B]$ the orbit of $B$ that is $[B]:=\{\{x^h,y^h\}\mid h\in H\}=\{B^h\mid h\in H\}$.

Given a multiset $M$ of elements of $G$ and $x\in G$, let $M_x:=[y\in M\mid y=x]$ be the sub-multiset of $M$ consisting of all copies of $x$. Write $w(x,M):=|M_x|.$

Now we consider an abelian group $(G,+)$, two subsets $A_1, A_2$ of $G$ and an automorphism group $H$ of $G$ such that
\begin{itemize}
\item $-Id\in H$ ($g^{-Id}:=-g$ for $g\in G$);
\item $H$ acts on $G\setminus A_1$ and on $G\setminus A_2$ respectively (i.e. $A_1$ and $A_2$ are union of orbits).
\end{itemize}
$H$ can be seen as an action on $G$.

\begin{Lemma}\label{wellDef1}
Let $B$ be a pair of elements of $G$ and $\mathcal{U}(B):=\bigcup_{\tilde{B}\in [B]} \tilde{B}$ be a  multiset. Then given any orbit $[z]\in G/H$, $w(z_1,\mathcal{U}(B))$ does not depend on the choice of $z_1\in [z]$.
\end{Lemma}

\proof
Let $z_1,z_2\in [z]$. Then there exists an $h\in H$ such that $z_2= z_1^h$ and we have
$$\mathcal{U}(B)=\bigcup_{\tilde{B}\in [B]} \tilde{B}=\bigcup_{\tilde{B}\in [B]} \tilde{B}^h=\mathcal{U}(B)^h.$$
Therefore $h$ induces a bijection between $\mathcal{U}(B)_{z_1}$ and $\mathcal{U}(B)_{z_2}$.
\endproof

\begin{Lemma}\label{wellDef2}
Let $B=\{x,y\}$ be a pair of elements of $G$ and $D(B):=\{x+y,x-y\}$. Let $\mathcal{D}(B):=\bigcup_{\tilde{B}\in [B]} D(\tilde{B})$ be a multiset. Then given any orbit $[z]\in G/H$, $w(z_1,\mathcal{D}(B))$ does not depend on the choice of $z_1\in [z]$.
\end{Lemma}

\proof
Let $z_1,z_2\in [z]$. Then there exists an $h\in H$ such that $z_2=z_1^h$. Since $h$ is an automorphism of $(G,+)$, we have
$$\mathcal{D}(B)=\bigcup_{\tilde{B}\in [B]} D(\tilde{B})=\bigcup_{\tilde{B}\in [B]} D(\tilde{B}^h)=\bigcup_{\tilde{B}\in [B]} D(\tilde{B})^h=\mathcal{D}(B)^h.$$
Therefore $h$ induces a bijection between $\mathcal{D}(B)_{z_1}$ and $\mathcal{D}(B)_{z_2}$.
\endproof

Let $G/H=\{[x_1],[x_2],\dots,[x_n]\}$ and ${{G}\choose{2}}/H=\{[B_1],[B_2],\dots,[B_m]\}$. Define a $2n\times m$ matrix $M=(M_{i,j})$, where
$$M_{i,j}=
\begin{cases}
w(x_i,\mathcal{U}(B_j)),\ \mbox{ if } 1 \leq i\leq n;\\
w(x_{i-n},\mathcal{D}(B_j)),\ \mbox{ if } n+1\leq i \leq 2n.
\end{cases}$$
The first $n$ rows of $M$ are labeled by $[x_i]$ for $1 \leq i\leq n$, and the last $n$ rows are labeled by $[x_{i-n}]$ for $n+1\leq i \leq 2n$. Define a column vector $J=(J_i)$ of length $2n$, where
$$J_{i}=
\begin{cases} 1,\ \mbox{ if } 1 \leq i\leq n,\  x_i\in G\setminus A_1;\\
1,\ \mbox{ if } n+1\leq i \leq 2n,\  x_{i-n}\in G\setminus A_2;\\
0,\ \mbox{ otherwise.}
\end{cases}$$

\begin{Proposition}\label{Kramer-Mesner}
Let $H$ be an automorphism group of $(G,+)$ such that $-Id \in H$ and let $A_1, A_2$ be subsets of $G$ which are union of orbits of the action of $G$ on $H$. If there exists a $0$-$1$ solution vector $X$ to $MX=J$, then there exists a PPS$(G,A_1,A_2)$.
\end{Proposition}

\proof
Since $MX=J$ has a $0$-$1$ solution vector $X=(X_j)$, if $[x_i]\subseteq G\setminus A_1$ for $1\leq i\leq n$, then there exists a unique orbit $[B_j]\in {{G}\choose{2}}/H$ such that $w(x_i,\mathcal{U}(B_j))=1$ and $X_j=1$. By Lemma \ref{wellDef1}, $w(x,\mathcal{U}(B_j))=1$ for any $x\in [x_i]$.

If $[x_{i-n}]\subseteq G\setminus A_2$ for $n+1\leq i\leq 2n$, then there exists a unique orbit $[B_j]\in {{G}\choose{2}}/H$ such that $w(x_i,\mathcal{D}(B_j))=1$ and $X_j=1$. By Lemma \ref{wellDef2}, $w(x,\mathcal{D}(B_j))=1$ for any $x\in [x_{i-n}]$.

If $[x_i]\subseteq A_1$ for $1\leq i\leq n$, then for any $1\leq j\leq m$ such that $X_j=1$, we have $w(x_i,\mathcal{U}(B_j))=0$. If $[x_{i-n}]\subseteq A_2$ for $n+1\leq i\leq 2n$, then for any $1\leq j\leq m$ such that $X_j=1$, we have $w(x_i,\mathcal{D}(B_j))=0$.

Write $${\cal S}_1=\bigcup\limits_{\substack{j=1\\ X_j=1}}^m[B_j].$$
It follows that $\bigcup_{B\in {\cal S}_1} B= G\setminus A_1$ and $\bigcup_{B\in {\cal S}_1} D(B)= G\setminus A_2.$

Since $-Id \in H$, for any $B\in {\cal S}_1$, we have $-B\in {\cal S}_1$. Now we define ${\cal S}$ to be a subfamily of ${\cal S}_1$ that contains exactly one element from each pair $\{B,-B\}$ where $B$ varies in $S_1$. Then ${\cal S}$ is a PPS$(G,A_1,A_2)$. \endproof


As an application of Proposition \ref{Kramer-Mesner} we obtain the following partitionable sets.

\begin{Theorem}\label{NewPS}
Let $q\equiv 3\pmod{4}$ be a prime and $3<q<200$. Let $p\equiv 3\pmod{4}$ be a prime, $p>3$ and $(p-1)\mid (q-1)$. Then there exists a PS$(pq)$, i.e., a ZCPS-Wh$(pq)$.
\end{Theorem}

\proof By Proposition \ref{prop:nece PS}, a PS$(pq)$, i.e., a ZCPS-Wh$(pq)$, exists only if $pq\equiv 1,5 \pmod{12}$. Thus it is required in the assumption that $p$ and $q$ are both greater than $3$. When $p=q$, the conclusion follows from Theorem \ref{prop:PS known}(2). For all the other values of $p$ and $q$, we list them in the following table.
\begin{center}
\begin{tabular}{|l|l|l|}
\hline $p$ & $q$ \\
\hline $7$ & $19,31,43,67,79,103,127,139,151,163,199$  \\
\hline $11$ & $31,71,131,151,191$  \\
\hline $19$ & $127,163,199$ \\
\hline $23$ & $67,199$ \\
\hline $31$ & $151$  \\
\hline $43$ & $127$ \\
\hline $47$ & $139$ \\
\hline $67$ & $199$ \\
\hline
\end{tabular}
\end{center}
To apply Proposition \ref{Kramer-Mesner} to find PS$(pq)$s, we take $G=\mathbb{Z}_{pq}$. Since $p$ and $q$ are coprime, $G$ is isomorphic to $\mathbb{Z}_{p}\times \mathbb{Z}_q$ under the mapping $\rho:x\mapsto (x\pmod{p},x\pmod{q})$, where $x\in G$. Let $\xi_{p}$ be an element of order $p-1$ in $\mathbb{Z}^*_{p}$ and $\xi_{q}$ be an element of order $p-1$ in $\mathbb{Z}^*_{q}$. Then $(\xi_{p},\xi_{q})$ is an element of order $p-1$ in $\mathbb{Z}^*_{p}\times \mathbb{Z}^*_q$. Clearly $(-1,-1)\in\langle(\xi_{p},\xi_{q})\rangle$, which is the multiplicative group generated by $(\xi_{p},\xi_{q})$ in $\mathbb{Z}^*_{p}\times \mathbb{Z}^*_q$. Take the element $\xi_{p,q}$ from $G$ such that $\xi_{p,q}$ corresponds to $(\xi_{p},\xi_{q})$ under the mapping $\rho$. Let $H=\langle\xi_{p,q}\rangle$. Then $-1\in H$ and $H$ is an automorphism group of $G$. By Proposition \ref{Kramer-Mesner}, it suffices to solve the system $MX=J$ by computer search.

For example, when $p=7$ and $q=19$, we take $\xi_7=3$ and $\xi_{19}=8$. Then $\xi_{7,19}=122$. Take $H=\langle122\rangle$. Then solving $MX=J$, we find 11 initial pairs:
$ \{ 122, 2 \}, \{ 100, 4 \}, \{ 78, 13 \}, \{ 67, 27 \}$, $\{ 56, 18 \}, \{ 46, 40 \}, \{ 79, 30 \}, \{ 112, 57 \}, \{ 101, 53 \},
  \{ 90, 63 \}, \{ 45, 34 \}.
 $
By the action of $H$, we obtain the following PS$(133)$:
\begin{center}\tabcolsep 0.022in
\begin{tabular}{lllllllll}
\underline{$\{122, 2\}$},& $\{121, 111\}$,& $\{132, 109\}$,& \underline{$\{100, 4\}$},& $\{97, 89\}$,& $\{130, 85\}$,& \underline{$\{78, 13\}$},& $\{73, 123\}$,& $\{128, 110\}$,\\
\underline{$\{67, 27\},$} & $ \{61, 102\}$,& $\{127, 75\}$,& \underline{$\{56, 18\}$}, & $\{49, 68\}$,& $\{126, 50\}$,& \underline{$\{46, 40\}$},& $\{26, 92\}$,& $\{113, 52\}$,\\
\underline{$\{79, 30\}$},& $\{62, 69\}$,& $\{116, 39\},$&  \underline{$ \{112, 57\}$},& $\{98, 38\}$,& $\{119, 114\}$,& \underline{$\{101, 53\}$},& $\{86, 82\}$, &$\{118, 29\}$,\\
\underline{$\{90, 63\}$}, &$\{74, 105\}$,& $\{117, 42\}$,& \underline{$\{45, 34\}, $}&
$ \{37, 25\}$, &$\{125, 124\}$.\\
\end{tabular}
\end{center}
The interested reader can get a copy of the solutions for all the other values of $p$ and $q$ from \cite{link}. \endproof

Proposition \ref{Kramer-Mesner} can be also used to look for almost partitionable sets. Here we list a few $APS$s which are obtained by using nontrivial group actions. Examples which use the action of $\{Id,-Id\}$ will be studied in the next paragraph.

\begin{Theorem}\label{NewAPS}
There is an APS$(v,v/3,v/3)$, i.e., a ZCPS-Wh$(v+1)$, for $v\in \{651,2343,4323\}$.
\end{Theorem}

\proof Let $v=3pq$, where $(v,p,q)\in\{(651,7,31),(2343,11,71),(4323,11,131)\}$. To apply Proposition \ref{Kramer-Mesner} to find APS$(3pq)$s, we take $G=\mathbb{Z}_{3pq}$. Note that $G$ is isomorphic to $\mathbb{Z}_{3}\times \mathbb{Z}_{p}\times \mathbb{Z}_q$ under the mapping $\rho:x\mapsto (x\pmod{3},x\pmod{p},x\pmod{q})$, where $x\in G$. Let $\xi_{p}$ be an element of order $p-1$ in $\mathbb{Z}^*_{p}$ and $\xi_{q}$ be an element of order $p-1$ in $\mathbb{Z}^*_{q}$. Then $(-1,\xi_{p},\xi_{q})$ is an element of order $p-1$ in $\mathbb{Z}^*_{3}\times \mathbb{Z}^*_{p}\times \mathbb{Z}^*_q$. Since $(p-1)/2$ is odd, we have $(-1,-1,-1)\in\langle(-1,\xi_{p},\xi_{q})\rangle$. Take the element $\xi_{p,q}$ from $G$ such that $\xi_{p,q}$ corresponds to $(-1,\xi_{p},\xi_{q})$ under the mapping $\rho$. Let $H=\langle\xi_{p,q}\rangle$. Then $-1\in H$ and $H$ is an automorphism group of $G$. By Proposition \ref{Kramer-Mesner}, it suffices to solve the system $MX=J$ by computer search.

For example, when $p=7$ and $q=31$, we take $\xi_7=5$ and $\xi_{31}=6$. Then $\xi_{7,31}=68$. Take $H=\langle68\rangle$. Solving $MX=J$, we find 54 initial pairs which are marked by underlines below.
By the action of $H$, we obtain the following APS$(651,217,217)$:
\begin{center}\tabcolsep 0.03in
\begin{tabular}{lllllllll}
\underline{$\{ 68, 3 \}$}, &$\{ 67, 204 \}$, &$\{ 650, 201 \}$,&
\underline{$\{ 136, 6 \}$}, &$\{ 134, 408 \}$, &$\{ 649, 402 \}$,&
\underline{$\{ 272, 9 \}$}, &$\{ 268, 612 \},$\\
$\{ 647, 603 \}$, &\underline{$\{ 340, 10 \}$}, &$\{ 335, 29 \}$, &$\{ 646, 19 \}$, &\underline{$\{ 476, 61 \}$}, &$\{ 469, 242 \}$, &$\{ 644, 181 \}$, &\underline{$\{ 544, 11 \},$}\\
$\{ 536, 97 \}$, &$\{ 643, 86 \}$, &\underline{$\{ 165, 24 \}$}, &$\{ 153, 330 \}$, &$\{ 639, 306 \}$,& \underline{$\{ 233, 178 \}$}, &$\{ 220, 386 \}$, &$\{ 638, 208 \},$\\
\underline{$
  \{ 301, 60 \}$},& $\{ 287, 174 \}$,& $\{ 637, 114 \}$,& \underline{$\{ 369, 127 \}$},& $\{ 354, 173 \}$,& $\{ 636, 46 \}$,& \underline{$\{ 437, 102 \}$},& $\{ 421, 426 \}$,\\

  $\{ 635, 324 \}$,& \underline{$\{ 505, 51 \}$}, & $\{ 488, 213 \}$,& $\{ 634, 162 \}$,& \underline{$\{ 58, 43 \}$}, & $\{ 38, 320 \}$,& $\{ 631, 277 \}$,& \underline{$\{ 573, 210 \}$},\\

  $\{ 555, 609 \}$,& $\{ 633, 399 \}$,& \underline{$\{ 126, 45 \}$},& $\{ 105, 456 \}$,& $\{ 630, 411 \}$,& \underline{$\{ 194, 89 \}$},& $\{ 172, 193 \}$,& $\{ 629, 104 \}$,&\\
  \underline{$\{ 398, 31 \}$},& $\{ 373, 155 \}$,& $\{ 626, 124 \}$,& \underline{$\{ 262, 185 \}$},& $\{ 239, 211 \}$,& $\{ 628, 26 \}$,& \underline{$\{ 534, 66 \}$},& $\{ 507, 582 \}$,&\\
  $\{ 624, 516 \}$,& \underline{$\{ 602, 142 \}$},& $\{ 574, 542 \}$,& $\{ 623, 400 \}$,& \underline{$\{ 87, 74 \}$},& $\{ 57, 475 \}$,& $\{ 621, 401 \}$,& \underline{$\{ 223, 92 \}$},&\\
  $\{ 191, 397 \}$,& $\{ 619, 305 \}$,& \underline{$\{ 291, 84 \}$},& $\{ 258, 504 \}$,& $\{ 618, 420 \}$,& \underline{$\{ 427, 184 \}$},& $\{ 392, 143 \}$,& $\{ 616, 610 \}$,&\\
  \underline{$\{ 563, 232 \}$},& $\{ 526, 152 \}$,& $\{ 614, 571 \}$,& \underline{$\{ 359, 55 \}$},& $\{ 325, 485 \}$,& $\{ 617, 430 \}$,& \underline{$\{ 145, 280 \}$},& $\{ 95, 161 \}$,&\\
  $\{ 601, 532 \}$,& \underline{$\{ 116, 183 \}$},& $\{ 76, 75 \}$,& $\{ 611, 543 \}$,& \underline{$\{ 495, 149 \}$},& $\{ 459, 367 \}$,& $\{ 615, 218 \}$,& \underline{$\{ 388, 79 \}$},&\\
  $\{ 344, 164 \}$,& $\{ 607, 85 \}$,& \underline{$\{ 349, 54 \}$},& $\{ 296, 417 \}$,& $\{ 598, 363 \}$,& \underline{$\{ 281, 285 \}$},& $\{ 229, 501 \}$,& $\{ 599, 216 \}$,&\\
  \underline{$\{ 592, 133 \}$},& $\{ 545, 581 \}$,& $\{ 604, 448 \}$,& \underline{$\{ 553, 72 \}$},& $\{ 497, 339 \}$,& $\{ 595, 267 \}$,& \underline{$\{ 271, 215 \}$},& $\{ 200, 298 \}$,&\\
  $\{ 580, 83 \}$,& \underline{$\{ 494, 177 \}$},& $\{ 391, 318 \}$,& $\{ 548, 141 \}$,& \underline{$\{ 513, 160 \}$},& $\{ 381, 464 \}$,& $\{ 519, 304 \}$,& \underline{$\{ 378, 308 \}$},&\\
  $\{ 315, 112 \}$,& $\{ 588, 455 \}$,& \underline{$\{ 484, 206 \}$},& $\{ 362, 337 \}$,& $\{ 529, 131 \}$,& \underline{$\{ 446, 207 \}$},& $\{ 382, 405 \}$,& $\{ 587, 198 \}$,&\\
  \underline{$\{ 357, 179 \}$},& $\{ 189, 454 \}$,& $\{ 483, 275 \}$,& \underline{$\{ 300, 91 \}$},& $\{ 219, 329 \}$,& $\{ 570, 238 \}$,& \underline{$\{ 523, 295 \}$},& $\{ 410, 530 \}$,&\\
  $\{ 538, 235 \}$,& \underline{$\{ 465, 99 \}$},& $\{ 372, 222 \}$,& $\{ 558, 123 \}$,& \underline{$\{ 261, 227 \}$},& $\{ 171, 463 \}$,& $\{ 561, 236 \}$,& \underline{$\{ 338, 266 \}$},&\\
  $\{ 199, 511 \}$,& $\{ 512, 245 \}$,& \underline{$\{ 407, 255 \}$},& $\{ 334, 414 \}$,& $\{ 578, 159 \}$,& \underline{$\{ 358, 274 \}$},& $\{ 257, 404 \}$,& $\{ 550, 130 \}$,&\\
  \underline{$\{ 310, 111 \}$},& $\{ 248, 387 \}$,& $\{ 589, 276 \}$,& \underline{$\{ 514, 190 \}$},& $\{ 449, 551 \}$,& $\{ 586, 361 \}$,& \underline{$\{ 319, 228 \}$},& $\{ 209, 531 \}$,&\\
  \end{tabular}
\end{center}
\begin{center}\tabcolsep 0.03in
\begin{tabular}{lllllllll}
  $\{ 541, 303 \}$,& \underline{$\{ 309, 169 \}$},& $\{ 180, 425 \}$,& $\{ 522, 256 \}$,& \underline{$\{ 368, 158 \}$},& $\{ 286, 328 \}$,& $\{ 569, 170 \}$,& \underline{$\{ 533, 151 \}$},&\\
  $\{ 439, 503 \}$,& $\{ 557, 352 \}$.\\
\end{tabular}
\end{center}
The interested reader can get a copy of the solutions for $v\in \{2343,4323\}$ from \cite{link}. \endproof

\begin{Theorem}\label{thm:<300}
Let $v\equiv 3\pmod{4}$, $v<300$ and $\alpha\beta\neq 0$. Then there exists an APS$(v,\alpha,\beta)$ if and only if
\begin{eqnarray*}
2\alpha^2-\beta^2\equiv
\left\{
\begin{array}{ll}
    \frac{v}{3} \ ({\rm mod }\ v), & \hbox{ \rm{if} $v\equiv 3\ ({\rm mod }\ 12)$,}\\[0.4em]
	0\ ({\rm mod }\ v), & \hbox{ \rm{if} $v\equiv 7,11 \ ({\rm mod }\ 12)$.}
\end{array}
\right.
\end{eqnarray*}
\end{Theorem}

\proof Recall that Lemma \ref{Nec} is obtained by refining Lemma \ref{nec}. By Lemma \ref{Nec}(1), there is no APS$(v)$ for $v\in\{63, 99, 171, 207, 279\}$; By Lemma \ref{Nec}(2), there is no APS$(v)$ for $v\in\{15,87,159$, $195\}$; By Lemma \ref{Nec}(3), there is no APS$(v)$ for $v\in\{11,19,43,55,59,67,83,95,107,131,139,143$, $163,179,211,215,227,247,251,283,295\}$. The case of $v=3$ is trivial. Thus for $v\equiv 3\pmod{4}$ and $v<300$, it remains to examine $v\in\{7,23,27,31,35,39,47,51,71,75,79,91,103,111,115,119,123$, $127,135,147,151,155,167, 175,183,187,191,199,203,219,223,231,235,239,243,255,259,263,267$, $271,275,287,291,299\}$.

By Corollary \ref{cor:APS from Marco}, there is an APS$(v,\alpha,\beta)$ for any prime $v\in \{7,23,31,47,71,127,151,167,191$, $263,271\}$ and $2\alpha^2-\beta^2\equiv0\pmod{v}$. By choosing proper values of $p$ and $q$ in Corollary \ref{cor:APS from recur}, an APS$(v,\alpha,\beta)$ exists for the following values of $v=pq$ and $2\alpha^2-\beta^2\equiv0\pmod{v}$:
\begin{center}
\begin{tabular}{|c|c|c|c|c|c|c|c|c|c|c}
\hline $v$ & $35$ & $91$ & $115$  & $155$ & $203$ & $235$ & $259$ & $299$ \\
\hline $q$ & $5$ & $13$ & $5$  & $5$& $29$& $5$& $37$& $13$ \\
\hline $p$ & $7$ & $7$ & $23$   & $31$& $7$& $47$& $7$& $23$ \\
\hline
\end{tabular} .
\end{center}
By Corollary \ref{cor:nec}, if $v\in\{39,111,183\}$, then an APS$(v,\alpha,\beta)$ exists only if $\alpha,\beta\in\{v/3,2v/3\}$. By Theorem \ref{SmallAPS}, there is an APS$(v,v/3,v/3)$ for any $v\in\{39,111,183\}$.

For $v\in\{27,51,75,79,103,119,123,135,147,175,187,199,219,223,231,239,243,255,267,275$, $287,291\}$, we will use Proposition \ref{Kramer-Mesner} and the action of the trivial group $\{Id,-Id\}$ to deal with all desired values of $\alpha$ and $\beta$ up to multiplier. Here we only write explicitly an APS$(v,\alpha,\alpha)$ for $(v,\alpha)\in \{(243,18),(255,85),(275,110)\}$, i.e., a ZCPS-Wh$(v+1)$, which are unknown before by Theorem \ref{SmallAPS}.
\begin{center}\tabcolsep 0.02in
\begin{tabular}{l}
\hline{APS$(243,18,18)$:}\\
\hline $\{1,2\}$, $\{63,99\}$, $\{7,13\}$, $\{3,29\}$, $\{8,37\}$, $\{5,54\}$, $\{10,68\}$, $\{9,84\}$, $\{11,101\}$,
$\{12,73\}$, $\{14,111\}$,\\$\{17,33\}$,
$\{16,103\}$, $\{15,77\}$, $\{19,85\}$, $\{20,53\}$, $\{21,93\}$, $\{25,119\}$,
$\{26,30\}$, $\{27,38\}$, $\{32,56\}$,\\$\{36,96\}$, $\{28,67\}$,
$\{35,98\}$, $\{39,106\}$, $\{34,118\}$, $\{42,79\}$,
$\{51,86\}$, $\{40,87\}$, $\{57,100\}$, $\{61,80\}$,\\$\{52,109\}$, $\{55,83\}$, $\{58,89\}$,
$\{72,102\}$, $\{114,121\}$, $\{41,95\}$, $\{66,88\}$, $\{45,115\}$, $\{49,74\}$, $\{4,6\}$,\\
$\{50,90\}$, $\{64,108\}$, $\{112,117\}$,
$\{46,120\}$, $\{94,107\}$, $\{47,81\}$, $\{62,113\}$, $\{43,91\}$, $\{70,97\}$,\\$\{44,82\}$,
$\{71,92\}$, $\{59,105\}$, $\{48,65\}$, $\{60,75\}$,
$\{22,31\}$, $\{23,78\}$, $\{24,76\}$, $\{69,110\}$, $\{104,116\}$.\\
\hline{APS$(255,85,85)$:}\\
\hline $\{1,2\}$, $\{3,9\}$, $\{4,17\}$, $\{5,27\}$, $\{6,37\}$, $\{7,52\}$, $\{8,62\}$, $\{10,105\}$, $\{11,83\}$,
$\{12,119\}$, $\{14,55\}$, \\ $\{13,93\}$,$\{15,123\}$, $\{21,61\}$, $\{16,87\}$, $\{18,29\}$, $\{19,33\}$,
$\{20,76\}$, $\{22,71\}$, $\{23,127\}$, $\{24,112\}$,\\ $\{31,66\}$, $\{28,56\}$, $\{25,89\}$, $\{30,45\}$,
$\{26,109\}$, $\{32,110\}$, $\{39,47\}$, $\{34,94\}$, $\{35,122\}$, $\{41,96\}$,\\ $\{49,114\}$, $\{36,98\}$,
$\{40,90\}$, $\{57,99\}$, $\{50,126\}$, $\{78,104\}$, $\{103,108\}$, $\{69,120\}$, $\{51,58\}$, \\ $\{64,91\}$,
$\{75,113\}$, $\{42,60\}$, $\{86,88\}$, $\{65,101\}$, $\{68,97\}$, $\{72,106\}$, $\{79,118\}$, $\{48,116\}$,\\
$\{115,124\}$, $\{53,73\}$, $\{74,107\}$, $\{63,82\}$, $\{46,70\}$, $\{117,121\}$, $\{92,102\}$, $\{44,67\}$,
$\{43,80\}$, \\$\{77,125\}$, $\{54,100\}$, $\{59,84\}$, $\{81,111\}$, $\{38,95\}$.\\

\hline{APS$(275,110,110)$:}\\
\hline $\{1, 2\}$, $\{3, 25\}$, $\{4, 6\}$, $\{5, 13\}$, $\{7, 36\}$, $\{8, 52\}$, $\{9, 71\}$, $\{10, 131\}$, $\{11, 90\}$, $\{12, 79\}$, $\{14,
  112\}$,\\ $\{15, 48\}$, $\{16, 73\}$, $\{17, 111\}$, $\{18, 121\}$, $\{19,
 53\}$, $\{20, 34\}$, $\{21, 28\}$, $\{22, 109\}$, $\{23, 119\}$, \\$\{24, 75\}$, $\{26, 96\}$,$\{27, 100\}$, $\{29, 84\}$, $\{30, 89\}$, $\{31, 54\}$, $\{32, 43\}$, $\{33, 59\}$, $\{35, 39\}$,
  $\{37, 127\}$, \\$\{38, 135\}$, $\{40, 128\}$, $\{41, 97\}$, $\{42, 88\}$, $\{44, 80\}$, $\{45, 126\}$, $\{46, 124\}$, $\{47, 82\}$, $\{49, 118\}$, \\$\{50, 116\}$, $\{51, 101\}$, $\{55, 108\}$, $\{56,
  133\}$, $\{57, 102\}$, $\{58, 77\}$, $\{60, 98\}$, $\{61, 132\}$,
  $\{62, 99\}$, \\$\{63, 69\}$, $\{64, 105\}$, $\{65, 117\}$, $\{66, 91\}$, $\{67, 83\}$, $\{68, 92\}$, $\{70, 85\}$, $\{72, 120\}$,
  $\{74, 106\}$, \\ $\{76, 123\}$,$\{78, 136\}$, $\{81, 94\}$, $\{86, 125\}$, $\{87, 104\}$, $\{93, 114\}$, $\{95, 115\}$, $\{103, 130\}$,
  $\{107, 137\}$, \\$\{113, 122\}$, $\{129, 134\}$.\\
\hline
\end{tabular}
\end{center}
The interested reader can get a copy of all the other data from \cite{link}. \qed

\section{Applications to optical orthogonal codes}\label{OOCSection}

The target of this section is to construct optical orthogonal codes using partitionable sets and almost partitionable sets.

A $(v,k,1)$-{\em optical orthogonal code} (OOC) is defined as a set ${\cal B}=\{B_1,B_2,\ldots,B_s\}$ of $k$-subsets (called {\em codewords}) of $\mathbb{Z}_v$ whose list of differences
$$\Delta {\cal B}:=\bigcup_{i=1}^s \Delta B_i:=\bigcup_{i=1}^s \ [x-y\mid x,y\in{\cal B}_i, x\not=y]$$
does not contain repeated elements. The number of codewords of an OOC is called the {\em size} of the OOC. In practice, a code with a large size is required.
The size of a $(v,k,1)$-OOC is upper bounded by the Johnson bound (cf. \cite{Yin98})
$$J(v,k,1):=\left\lfloor\frac{v-1}{k(k-1)}\right\rfloor.$$
The {\em difference leave} or briefly {\em leave} $L(\cal{B})$ of a $(v,k,1)$-OOC $\cal{B}$ is defined by the set of missing differences $L(\cal{B}):=$ $\mathbb{Z}_v\setminus (\{0\}\cup\Delta {\cal B})$.



Strong difference families have been used to construct OOCs in \cite{b18,ccfw}. The idea of strong difference families was also implicitly used in \cite{MaChang05,yyl} to construct OOCs. Suppose $(G,+)$ is a finite group and let $\Sigma=[F_1,F_2,\dots,F_s]$ be a family of multisets of size $k$ of $G$, where $F_i=[f_{i,1},f_{i,2},\ldots,f_{i,k}]$ for $1\leq i\leq s$. $\Sigma$ is said to be a $(G,k,\mu)$ {\em strong difference family} if the list
$$\Delta \Sigma:=\bigcup_{i=1}^s \Delta F_i:=\bigcup_{i=1}^s\ [f_{i,a}-f_{i,b}: 1\leq a,b\leq k, a\not=b]=\underline{\mu} G,$$
i.e., every element of $G$ $($0 included$)$ appears exactly $\mu$ times in the multiset $\Delta \Sigma$. The members of $\Sigma$ are called {\em base blocks}.

The concept of strong difference family was introduced by Buratti in \cite{b99} and revisited in \cite{m}. Similar to what has been done in \cite{ccfw,cfw,cfw2}, here we will focus on three particular SDFs with some special ``patterns" and we will look for second components: the main ingredients for this purpose will be given by PSs and APSs. More precisely we will use the following three SDFs:
\begin{center}\tabcolsep 0.05in
\begin{tabular}{ll}
$(\mathbb{Z}_3,4,4)$-SDF: & $\Sigma_1=[[1,1,-1,-1]]$.\\
$(\mathbb{Z}_5,5,4)$-SDF: & $\Sigma_2=[[0,1,1,-1,-1]]$.\\
$(\mathbb{Z}_{45},5,4)$-SDF: & $\Sigma_3=[[0,1,1,-1,-1]$, $[0,3,7,13,30]$, $[0,3,7,13,30]$, $[0,3,7,13,30]$,\\
& $[0,3,7,13,30]$, $[0,5,14,26,34]$, $[0,5,14,26,34]$, $[0,5,14,26,34],$ $[0,5,14,26,34]]$.
\end{tabular}
\end{center}

\subsection{Families of OOCs achieving the Johnson bound}

\begin{Theorem}\label{OOC4}
Suppose there exists a PS$(v)$ or an APS$(v,\alpha,\beta)$. Then there exist
\begin{itemize}
\item[$1)$] a $(3v,4,1)$-OOC with $J(3v,4,1)=\lfloor(v-1)/4\rfloor$ codewords whenever $\gcd(v,6)=1$, and
\item[$2)$] a $(5v,5,1)$-OOC with $J(5v,5,1)=\lfloor(v-1)/4\rfloor$ codewords whenever $\gcd(v,10)=1$.
\end{itemize}
\end{Theorem}

\proof
Let $\mathcal{S}$ be a PS$(v)$ or an APS$(v,\alpha,\beta)$. We here construct a $(3v,4,1)$-OOC (resp. $(5v,5,1)$-OOC) on $\mathbb{Z}_3\times \mathbb{Z}_v\cong \mathbb{Z}_{3v}$ (resp. $\mathbb{Z}_5\times \mathbb{Z}_v\cong \mathbb{Z}_{5v}$). Let
$$\mathcal{F}_1=\{\{(1,x),(1,-x),(-1,y),(-1,-y)\}\mid \{x,y\}\in \mathcal{S}\}$$
and
$$\mathcal{F}_2=\{\{(0,0),(1,x),(1,-x),(-1,y),(-1,-y)\}\mid \{x,y\}\in \mathcal{S}\}.$$
Then $\mathcal{F}_1$ is a $(3v,4,1)$-OOC with $\lfloor(v-1)/4\rfloor$ codewords whenever $\gcd(v,6)=1$ and $\mathcal{F}_2$ is a $(5v,5,1)$-OOC with $\lfloor(v-1)/4\rfloor$ codewords whenever $\gcd(v,10)=1$.

It is readily checked that $\Delta\mathcal{F}_1=\bigcup_{i\in \mathbb{Z}_3} \{i\}\times D_i$
and $\Delta\mathcal{F}_2=\bigcup_{i\in \mathbb{Z}_5} \{i\}\times D'_i,$
where
\begin{itemize}
\item[] $D_0=D'_0= \bigcup_{\{x,y\}\in \mathcal{S}} \pm\{2x,2y\};$
\item[] $D_1=D_{-1}=D'_2=D'_{-2}=\bigcup_{\{x,y\}\in \mathcal{S}} \pm\{x-y,x+y\};$
\item[] $D'_1=D'_{-1}=\bigcup_{\{x,y\}\in \mathcal{S}} \pm\{x,y\}.$
\end{itemize}

If $\mathcal{S}$ is a PS$(v)$, then the leave of $\mathcal{F}_1$ is $L(\mathcal{F}_1)=(\mathbb{Z}_3\times \mathbb{Z}_v)\setminus \Delta \mathcal{F}_1=\mathbb{Z}_3\times\{0\}$, and so $\mathcal{F}_1$ forms a $(3v,4,1)$-OOC. Similarly, $L(\mathcal{F}_2)=(\mathbb{Z}_5\times \mathbb{Z}_v)\setminus \Delta \mathcal{F}_2=\mathbb{Z}_5\times\{0\}$, so $\mathcal{F}_2$ forms a $(5v,5,1)$-OOC.

If $\mathcal{S}$ is an APS$(v,\alpha,\beta)$, then $L(\mathcal{F}_1)= (\{0\}\times \{0,2\alpha,-2\alpha\}) \cup(\{1,-1\}\times\{0,\beta,-\beta\})$, which is of size 9, and so $\mathcal{F}_1$ is a $(3v,4,1)$-OOC. Similarly, $L(\mathcal{F}_2)= (\{0\}\times \{0,2\alpha,-2\alpha\}) \cup(\{1,-1\}\times\{0,\alpha,-\alpha\})\cup (\{2,-2\}\times\{0,\beta,-\beta\})$, which is of size 15, and so $\mathcal{F}_2$ is a $(5v,4,1)$-OOC. \endproof

We remark that in the proof of Theorem \ref{OOC4}, strong different families are employed implicitly since the first coordinates of elements in ${\cal F}_1$ and ${\cal F}_2$ form a $(\mathbb{Z}_3,4,4)$-SDF and a $(\mathbb{Z}_5,5,4)$-SDF, respectively. We also remark that if the input APS in Theorem \ref{OOC4} is an APS$(p,1,\theta-1)$ that is from Lemma \ref{Marco}, then the resulting OOCs are those of Theorem $3.1$ of \cite{b18}.

\begin{Theorem}\label{OOC45}
Let $v\equiv 1 \pmod{4}$ and $\gcd(v,45)=1$. If there exists a PS$(v)$, then there exists a  $(45v,5,1)$-OOC with $J(45v,5,1)=(9v-1)/4$ codewords.
\end{Theorem}

\proof
Let $\mathcal{S}$ be a PS$(v)$. We first construct the following family $\cal F$ of codewords on $\mathbb{Z}_{45}\times \mathbb{Z}_v\cong \mathbb{Z}_{45v}$, whose first coordinates form a $(\mathbb{Z}_{45},5,4)$-SDF:
\begin{center}
\begin{tabular}{lll}
$\mathcal{F}=$&$\{\{(0,0),(1,x),(1,-x),(-1,y),(-1,-y)\},$\\&
$ \{(0,0),(3,x),(7,2x),(13,3x),(30,4x)\}, $\\&
$ \{(0,0),(3,-x),(7,-2x),(13,-3x),(30,-4x)\},$\\&
$ \{(0,0),(3,y),(7,2y),(13,3y),(30,4y)\},$\\&
$ \{(0,0),(3,-y),(7,-2y),(13,-3y),(30,-4y)\},$\\&
$ \{(0,0),(5,x),(14,2x),(26,3x),(34,4x)\},$\\&
$ \{(0,0),(5,-x),(14,-2x),(26,-3x),(34,-4x)\},$\\&
$ \{(0,0),(5,y),(14,2y),(26,3y),(34,4y)\},$\\&
$ \{(0,0),(5,-y),(14,-2y),(26,-3y),(34,-4y)\}\mid \{x,y\}\in \mathcal{S}\}.$
\end{tabular}
\end{center}
It is readily checked that $\Delta\mathcal{F}=\bigcup_{i\in \mathbb{Z}_{45}} \{i\}\times D_i,$
where
\begin{itemize}
\item[] $D_0= \bigcup_{\{x,y\}\in \mathcal{S}} \pm\{2x,2y\};$
\item[] $D_1=D_{-1}=\bigcup_{\{x,y\}\in \mathcal{S}} \pm\{x,y\};$
\item[] $D_2=D_{-2}=\bigcup_{\{x,y\}\in \mathcal{S}} \pm\{x-y,x+y\};$
\item[] $D_i=D_{-i}=\bigcup_{\{x,y\}\in \mathcal{S}} \pm\{\delta_ix,\delta_iy\}$ for $3\leq i\leq 22$, where $\delta_i$ is an invertible element of $\mathbb{Z}_{v}$.
\end{itemize}
Since $\mathcal{S}$ is a PS$(v)$, then the leave of $\mathcal{F}$ is $L(\mathcal{F})=\mathbb{Z}_{45}\times\{0\}$. Furthermore, take
$$\mathcal{E}=\left\{\{(0,0),(1,0),(3,0),(29,0),(35,0)\},\{(0,0),(5,0),(20,0),(27,0),(41,0)\}\right\}.$$
Then $\mathcal{F}\cup\mathcal{E}$ forms a $(45v,5,1)$-OOC with $(9v-1)/4$ codewords.. \endproof

Applying Theorem \ref{OOC45} with a PS$(p)$ where $p\equiv 1\pmod{4}$ is a prime, we give another proof of Lemma 2.8 in \cite{MaChang04}.

\subsection{Families of OOCs achieving the Johnson bound minus one}

In this subsection, we shall provide some families of OOCs which fail to achieve the Johnson bound for just one codeword.

\begin{Lemma}\label{PPSinvertible}
Let $p\equiv 3\pmod{4}$ and $q\equiv 3 \pmod{4}$ be primes such that $p>q>3$. Then there exists a PPS$(\mathbb{Z}_{pq}, p\mathbb{Z}_{pq}\cup q\mathbb{Z}_{pq}, p\mathbb{Z}_{pq}\cup q\mathbb{Z}_{pq})$.
\end{Lemma}

\proof Since $p$ and $q$ are two distinct primes, $\mathbb{Z}_{pq}\cong \mathbb{Z}_{p}\times \mathbb{Z}_{q}$. Let $\Box_p$ and $\Box_q$ be the sets of all nonzero squares in $\mathbb{Z}_{p}$ and $\mathbb{Z}_{q}$, respectively. Let $\boxtimes_p$ and $\boxtimes_q$ be the sets of all non-squares in $\mathbb{Z}_{p}$ and $\mathbb{Z}_{q}$, respectively. Write
$$\mathcal{S}:=\left\{\{(s_1,s_2),(x_1s_1,x_2s_2)\} \mid s_1 \in \Box_p, s_2 \in \Box_q\right\},$$
where
\begin{center}
\begin{tabular}{ll}
$\begin{cases}
x_1\in \Box_p,\\
(1+x_1)(1-x_1)\in \Box_p,\\
\end{cases}$ and &
$\begin{cases}
x_2\in\boxtimes_q,\\
(1+x_2)(1-x_2)\in \boxtimes_q.
\end{cases}$
\end{tabular}
\end{center}
The required $x_1$ and $x_2$ exist for all primes $p\equiv 3\pmod{4}$ and $q\equiv 3 \pmod{4}$ that are greater than $3$ (note that $-1\in \boxtimes_p$ and $-1\in \boxtimes_q$). For example, we can take


\begin{eqnarray*}
(x_1,x_2)=
\left\{
\begin{array}{ll}
    (1/3,-2), & \hbox{ \rm{if} $2\in\Box_q$ \rm{and} $3\in \Box_q$,}\\[0.4em]
    (2,3), & \hbox{ \rm{if} $2\in\Box_q$ \rm{and} $3\in\boxtimes_q$,}\\[0.4em]
    (3,2), & \hbox{ \rm{if} $2\in\boxtimes_q$ \rm{and} $3\in \Box_q$,}\\[0.4em]
    (-3,1/2), & \hbox{ \rm{if} $2\in\boxtimes_q$ \rm{and} $3\in\boxtimes_q$.}
\end{array}
\right.
\end{eqnarray*}



\noindent It is readily checked that
\begin{align*}
& \bigcup_{\substack{s_1 \in \Box_p\\ s_2 \in \Box_q}}\pm\{(s_1,s_2),(x_1s_1,x_2s_2)\}
\\ = & \bigcup_{\substack{s_1 \in \Box_p\\ s_2 \in \Box_q}}\pm\left\{\left((1+x_1)s_1,(1+x_2)s_2\right),((1-x_1)s_1,(1-x_2)s_2)\right\}
\\ = &\ (\mathbb{Z}_p\times \mathbb{Z}_q)\setminus\left((\{0\}\times \mathbb{Z}_q)\cup(\mathbb{Z}_p\times \{0\})\right).
\end{align*}
Therefore, $\mathcal{S}$ forms a PPS$(\mathbb{Z}_p\times \mathbb{Z}_q,A,A)$, where $A= (\{0\}\times \mathbb{Z}_q)\cup(\mathbb{Z}_p\times \{0\})$. That is a PPS$(\mathbb{Z}_{pq}, p\mathbb{Z}_{pq}\cup q\mathbb{Z}_{pq}, p\mathbb{Z}_{pq}\cup q\mathbb{Z}_{pq})$.
\endproof

\begin{Lemma}\label{PPS2}
Suppose there exist an APS$(p,\alpha_1,\beta_1)$ and an APS$(q,\alpha_2,\beta_2)$, where $p$ and $q$ are primes such that $p>q>3$. Then there exists a PPS$(\mathbb{Z}_{pq},\{0,\pm q\alpha_1,\pm p\alpha_2\},\{0,\pm q\beta_1,\pm p\beta_2\})$.
\end{Lemma}

\proof Let $\mathcal{S}$ be a PPS$(\mathbb{Z}_{pq}, p\mathbb{Z}_{pq}\cup q\mathbb{Z}_{pq}, p\mathbb{Z}_{pq}\cup q\mathbb{Z}_{pq})$, which exists by Lemma \ref{PPSinvertible}. Let $\mathcal{S}_p$ be an APS$(p,\alpha_1,\beta_1)$ and write $q\mathcal{S}_p:=\{\{qx,qy\}\mid \{x,y\}\in\mathcal{S}_p\}$. Then
\begin{enumerate}
\item[] $\bigcup_{\{x,y\}\in \mathcal{S}_p}\pm \{qx, qy\}=q\mathbb{Z}_{pq}\setminus\{0, \pm q\alpha_1\}$;
\item[] $\bigcup_{\{x,y\}\in \mathcal{S}_p}\pm \{qx-qy, qx+qy\}=q\mathbb{Z}_{pq}\setminus\{0, \pm q\beta_1\}$.
\end{enumerate}
Similarly, let $\mathcal{S}_q$ be an APS$(q,\alpha_2,\beta_2)$ and write $p\mathcal{S}_q:=\{\{px,py\}\mid \{x,y\}\in\mathcal{S}_q\}$. Then
\begin{enumerate}
\item[] $\bigcup_{\{x,y\}\in \mathcal{S}_q}\pm \{px, py\}=p\mathbb{Z}_{pq}\setminus\{0, \pm p\alpha_2\}$;
\item[] $\bigcup_{\{x,y\}\in \mathcal{S}_q}\pm \{px-py, px+py\}=p\mathbb{Z}_{pq}\setminus\{0, \pm p\beta_2\}$.
\end{enumerate}
Then $\mathcal{S}\cup q\mathcal{S}_p\cup p\mathcal{S}_q$ forms a PPS$(\mathbb{Z}_{pq},\{0,\pm q\alpha_1,\pm p\alpha_2\},\{0,\pm q\beta_1,\pm p\beta_2\})$. \endproof

\begin{Theorem}\label{MaxOOC1}
Suppose there exist an APS$(p,\alpha_1,\beta_1)$ and an APS$(q,\alpha_2,\beta_2)$, where $p$ and $q$ are primes such that $p>q>3$. Then there exist
\begin{itemize}
\item[$1)$] a $(3pq,4,1)$-OOC $\mathcal{F}_1$ with $J(3pq,4,1)-1=(pq-5)/4$ codewords, and
\item[$2)$] a $(5pq,5,1)$-OOC $\mathcal{F}_2$ with $J(5pq,5,1)-1=(pq-5)/4$ codewords.
\end{itemize}
Furthermore, there is no other $(3pq,4,1)$-OOC $($resp. $(5pq,5,1)$-OOC$)$ $\mathcal{A}$ such that $\mathcal{A}$ properly contains $\mathcal{F}_1$ $($resp. $\mathcal{F}_2$$)$.
\end{Theorem}

\proof We here construct a $(3pq,4,1)$-OOC (resp. $(5pq,5,1)$-OOC) on $\mathbb{Z}_3\times \mathbb{Z}_{pq}\cong \mathbb{Z}_{3pq}$ (resp. $\mathbb{Z}_5\times \mathbb{Z}_{pq}\cong \mathbb{Z}_{5pq}$). Let $\mathcal{S}$ be a PPS$(\mathbb{Z}_{pq},\{0,\pm q\alpha_1,\pm p\alpha_2\},\{0,\pm q\beta_1,\pm p\beta_2\})$, which exists by Lemma \ref{PPS2}. Let
$$\mathcal{F}_1=\{\{(1,x),(1,-x),(-1,y),(-1,-y)\}\mid \{x,y\}\in \mathcal{S}\}$$
and
$$\mathcal{F}_2=\{\{(0,0),(1,x),(1,-x),(-1,y),(-1,-y)\}\mid \{x,y\}\in \mathcal{S}\}.$$
Then
$\Delta\mathcal{F}_1=\bigcup_{i\in \mathbb{Z}_3} \{i\}\times D_i$
and
$\Delta\mathcal{F}_2=\bigcup_{i\in \mathbb{Z}_5} \{i\}\times D'_i,$
where
\begin{itemize}
\item[] $D_0=D'_0=\mathbb{Z}_{pq}\setminus \{0,\pm 2q\alpha_1,\pm 2p\alpha_2\}$;
\item[] $D_{-1}=D_{1}=D'_{-2}=D'_{2}=\mathbb{Z}_{pq}\setminus \{0,\pm q\beta_1,\pm p\beta_2\}$;
\item[] $D'_{-1}=D'_{1}=\mathbb{Z}_{pq}\setminus \{0,\pm q\alpha_1,\pm p\alpha_2\}$.
\end{itemize}
The leave of $\mathcal{F}_1$ is
$$L(\mathcal{F}_1)=(\{0\}\times\{0,\pm 2q\alpha_1,\pm 2p\alpha_2\})\cup (\{1,-1\}\times \{0,\pm q\beta_1,\pm p\beta_2\})$$
and the leave of $\mathcal{F}_2$ is
$$L(\mathcal{F}_2)=(\{0\}\times\{0,\pm 2q\alpha_1,\pm 2p\alpha_2\})\cup (\{1,-1\}\times \{0,\pm q\alpha_1,\pm p\alpha_2\})\cup (\{2,-2\}\times \{0,\pm q\beta_1,\pm p\beta_2\}).$$
Therefore $|L(\mathcal{F}_1)|=15$ which yields that $\mathcal{F}_1$ is a $(3pq,4,1)$-OOC failing to reach the Johnson bound for exactly one codeword. Similarly, $|L(\mathcal{F}_2)|=25$ which yields that $\mathcal{F}_2$ is a $(5pq,5,1)$-OOC failing to reach the Johnson bound for exactly one codeword.

It remains to show that $\mathcal{F}_1$ and $\mathcal{F}_2$ cannot be extended.
Assume that $\mathcal{F}_1$ could be extended by adding a new codeword $B=\{v_1,v_2,v_3,v_4\}$ on $\mathbb{Z}_3\times \mathbb{Z}_{pq}$ satisfying
$$\Delta B\subseteq L(\mathcal{F}_1)\setminus\{(0,0)\}=(\{0\}\times\{\pm 2q\alpha_1,\pm 2p\alpha_2\})\cup ( \{1,-1\}\times \{0,\pm q\beta_1,\pm p\beta_2\}).$$
Note that each element in $\{\pm 2q\alpha_1,\pm q\beta_1\}$ is not divisible by $p$, and each element in $\{\pm 2p\alpha_2, \pm p\beta_2\}$ is not divisible by $q$. It follows that the second coordinates of $v_1,v_2,v_3,v_4$ must be divisible by $p$ (or by $q$) at the same time; otherwise, there would be some difference $v_i-v_j\not\in L(\mathcal{F}_1)\setminus\{(0,0)\}$. Hence $\Delta B$ contributes $12$ differences whose second coordinates are divisible by $p$ (or by $q$) at the same time, but this is absurd because $L(\mathcal{F}_1)\setminus\{(0,0)\}$ contains only $8$ such differences. Therefore $\mathcal{F}_1$ cannot be extended.

Similarly, assume that $\mathcal{F}_2$ could be extended by adding a new codeword $B=\{v_1,v_2,v_3,v_4,v_5\}$ on $\mathbb{Z}_5\times \mathbb{Z}_{pq}$ satisfying
$$\Delta B\subseteq (\{0\}\times\{\pm 2q\alpha_1,\pm 2p\alpha_2\})\cup (\{1,-1\}\times \{0,\pm q\alpha_1,\pm p\alpha_2\})\cup (\{2,-2\}\times \{0,\pm q\beta_1,\pm p\beta_2\}).$$
Note that each element in $\{\pm 2q\alpha_1,\pm q\alpha_1,\pm q\beta_1\}$ is not divisible by $p$, and each element in $\{\pm 2p\alpha_2,$ $\pm p\alpha_2, \pm p\beta_2\}$ is not divisible by $q$. It follows that the second coordinates of $v_1,v_2,v_3,v_4,v_5$ must be divisible by $p$ (or by $q$) at the same time, and so $\Delta B$ contributes $20$ differences whose second coordinates are divisible by $p$ (or by $q$) at the same time. A contradiction occurs because $L(\mathcal{F}_2)\setminus\{(0,0)\}$ contains only $14$ such differences. Therefore $\mathcal{F}_2$ cannot be extended. \endproof

Let $U(p^2)$ be the group of units in $\mathbb Z/p^2$. When $p\equiv 7 \pmod{8}$ is a prime, $2$ is a nonzero square in $U(p^2)$. In what follows we use $\sqrt{2}$ to represent one of the square roots of 2 in  $U(p^2)$.

\begin{Lemma}\label{SilverPPS}
Let $p\equiv 7 \pmod{8}$ be a prime. If $\theta=1+\sqrt{2}$ generates $U(p^2)/\{1,-1\}$, then there exists a PPS$(\mathbb{Z}_{p^2},\{0,\pm 1,\pm p\},\{0,\pm (\theta-1),\pm p(\theta-1)\})$.
\end{Lemma}

\proof Consider the set of pairs $\mathcal{S}_1=\{\{\theta^{2i-1},\theta^{2i}\}\mid 1\leq i\leq(p^2-p-2)/4\}$. Since $\theta$ generates $U(p^2)/\{1,-1\}$ and satisfies the equation $x(x-1)=x+1$, we have
$$\bigcup_{i=1}^{(p^2-p-2)/4}\pm\{\theta^{2i-1},\theta^{2i}\}=U(p^2)\setminus \{\pm 1\},$$
and
$$\bigcup_{i=1}^{(p^2-p-2)/4}\pm\{\theta^{2i-1}(\theta-1),\theta^{2i-1}(\theta+1)\}=U(p^2)\setminus\{\pm (\theta-1)\}.$$
Since the powers of $\theta$ cover, up to the sign, all the classes modulo $p^2$ that are coprime with $p$, those powers also cover (up to the sign) all the nonzero classes modulo $p$. Thus $\theta$ generates $\mathbb{Z}_{p}^*/\{1,-1\}$. Let $\mathcal{S}_2=\{\{p\theta^{2i-1},p\theta^{2i}\}\mid 1\leq i\leq(p-3)/4\}$. We have
$$\bigcup_{i=1}^{(p-3)/4}\pm\{p\theta^{2i-1},p\theta^{2i}\}=p\mathbb{Z}_{p^2}\setminus\{0,\pm p\},$$
and
$$\bigcup_{i=1}^{(p-3)/4}\pm\{p\theta^{2i-1}(\theta-1),p\theta^{2i-1}(\theta+1)\}=p\mathbb{Z}_{p^2}\setminus\{0,\pm p(\theta-1)\}.$$
Therefore, $\mathcal{S}_1\cup \mathcal{S}_2$ forms a PPS$(\mathbb{Z}_{p^2},\{0,\pm 1,\pm p\},\{0,\pm (\theta-1),\pm p(\theta-1)\})$.
\endproof

\begin{Lemma}\label{necsquare}
Let $p\equiv 7 \pmod{8}$ be a prime. If there is a PPS$(\mathbb{Z}_{p^2},\{0,\pm \alpha,\pm p\alpha\},\{0,\pm \beta,\pm p\beta\})$, then
$2\alpha^2-\beta^2\equiv 0\ ({\rm mod }\ p^2)$.
\end{Lemma}

\proof Let $\mathcal{S}$ be a PPS$(\mathbb{Z}_{p^2},\{0,\pm \alpha,\pm p\alpha\},\{0,\pm \beta,\pm p\beta\})$. By the definition of PPS, we can compute the sum of squares of all elements in $\mathbb{Z}_{p^2}$ in three ways as follows:
$$2\alpha^2+2p^2\alpha^2+2\sum_{\{x,y\}\in \mathcal{S}}(x^2+y^2)\equiv 2\beta^2+2p^2\beta^2+4\sum_{\{x,y\}\in \mathcal{S}}(x^2+y^2)\equiv \sum_{j=1}^{p^2-1}j^2\ ({\rm mod}\ p^2).$$
It follows that
\begin{align*}
&4\alpha^2+4p^2\alpha^2-2\beta^2-2p^2\beta^2 \\ &=2\left(2\alpha^2+2p^2\alpha^2+2\sum\limits_{\{x,y\}\in S}(x^2+y^2)\right)-\left(2\beta^2+2p^2\beta^2+4\sum\limits_{\{x,y\}\in S}(x^2+y^2)\right)
\\ &\equiv  2\sum\limits_{j=1}^{p^2-1}j^2-\sum\limits_{j=1}^{p^2-1}j^2  \equiv \sum\limits_{j=1}^{p^2-1}j^2 \equiv \frac{(p^2-1)p^2(2p^2-1)}{6}\equiv 0 \ ({\rm mod}\ p^2).
\end{align*}
Hence $2\alpha^2-\beta^2\equiv 0\ ({\rm mod }\ p^2)$.
\endproof

\begin{Theorem}\label{thm:PPS}
Suppose that $p\equiv 7 \pmod{8}$ is a prime, and $\theta=1+\sqrt{2}$ generates $U(p^2)/\{1,-1\}$. Then a PPS$(\mathbb{Z}_{p^2},\{0,\pm \alpha,\pm p\alpha\},\{0,\pm \beta,\pm p\beta\})$ exists if and only if $2\alpha^2-\beta^2\equiv0\pmod{p^2}$.
\end{Theorem}

\proof The necessity follows from Lemma \ref{necsquare}. It remains to examine the sufficiency. Assume that $\beta \equiv \pm \sqrt{2}\alpha \pmod{p^2}$. Let $\cal S$ be a PPS$(\mathbb{Z}_{p^2},\{0,\pm 1,\pm p\},\{0,\pm (\theta-1),\pm p(\theta-1)\})$ given by Lemma \ref{SilverPPS}. Then
$\{\{\alpha x,\alpha y\} \mid\{x,y\}\in{\cal S}\}$
is a PPS$(\mathbb{Z}_{p^2},\{0,\pm \alpha,\pm p\alpha\},\{0,\pm \beta,\pm p\beta\})$. \endproof

\begin{Theorem}\label{MaxOOC2}
Let $p\equiv 7 \pmod{8}$ be a prime. If $\theta=1+\sqrt{2}$ generates $U(p^2)/\{1,-1\}$, then there exist
\begin{itemize}
\item[$1)$] a $(3p^2,4,1)$-OOC $\mathcal{F}_1$ with $J(3p^2,4,1)-1=(p^2-5)/4$ codewords, and
\item[$2)$] a $(5p^2,5,1)$-OOC $\mathcal{F}_2$ with $J(5p^2,5,1)-1=(p^2-5)/4$ codewords.
\end{itemize}
Furthermore, there is no other $(3p^2,4,1)$-OOC $($resp. $(5p^2,5,1)$-OOC$)$ $\mathcal{A}$ such that $\mathcal{A}$ properly contains $\mathcal{F}_1$ $($resp. $\mathcal{F}_2$$)$.
\end{Theorem}

\proof We here construct a $(3p^2,4,1)$-OOC (resp. $(5p^2,5,1)$-OOC) on $\mathbb{Z}_3\times \mathbb{Z}_{p^2}\cong \mathbb{Z}_{3p^2}$ (resp. $\mathbb{Z}_5\times \mathbb{Z}_{p^2}\cong \mathbb{Z}_{5p^2}$). Let $\mathcal{S}$ be a PPS$(\mathbb{Z}_{p^2},\{0,\pm 1,\pm p\},\{0,\pm (\theta-1),\pm p(\theta-1)\})$ given by Lemma \ref{SilverPPS}. Let
$$\mathcal{F}_1=\{\{(1,x),(1,-x),(-1,y),(-1,-y)\}\mid \{x,y\}\in \mathcal{S}\}$$
and
$$\mathcal{F}_2=\{\{(0,0),(1,x),(1,-x),(-1,y),(-1,-y)\}\mid \{x,y\}\in \mathcal{S}\}.$$
Then $\Delta\mathcal{F}_1=\bigcup_{i\in \mathbb{Z}_3} \{i\}\times D_i$
and
$\Delta\mathcal{F}_2=\bigcup_{i\in \mathbb{Z}_5} \{i\}\times D'_i,$
where
\begin{itemize}
\item[] $D_0=D'_0=\mathbb{Z}_{p^2}\setminus \{0,\pm 2,\pm 2p\}$;
\item[] $D_{-1}=D_{1}=D'_{-2}=D'_{2}=\mathbb{Z}_{p^2}\setminus\{0,\pm (\theta-1),\pm p(\theta-1)\}$;
\item[] $D'_{-1}=D'_{1}=\mathbb{Z}_{p^2}\setminus\{0,\pm 1,\pm p\}$.
\end{itemize}
The leave of $\mathcal{F}_1$ is
$$L(\mathcal{F}_1)=(\{0\}\times\{0,\pm 2,\pm 2p\})\cup (\{1,-1\}\times \{0,\pm \sqrt{2},\pm p\sqrt{2}\})$$
and the leave of $\mathcal{F}_2$ is
$$L(\mathcal{F}_2)=(\{0\}\times\{0,\pm 2,\pm 2p\})\cup (\{1,-1\}\times \{0,\pm 1,\pm p\})\cup (\{2,-2\}\times \{0,\pm \sqrt{2},\pm p\sqrt{2}\}).$$
Therefore $|L(\mathcal{F}_1)|=15$ which yields that $\mathcal{F}_1$ is a $(3p^2,4,1)$-OOC failing to reach the Johnson bound for exactly one codeword. Similarly, $|L(\mathcal{F}_2)|=25$ which yields that $\mathcal{F}_2$ is a $(5p^2,5,1)$-OOC failing to reach the Johnson bound for exactly one codeword.

It remains to show that $\mathcal{F}_1$ and $\mathcal{F}_2$ cannot be extended.
Assume that $\mathcal{F}_1$ could be extended by adding a new codeword $B=\{v_1,v_2,v_3,v_4\}$ on $\mathbb{Z}_3\times \mathbb{Z}_{p^2}$ satisfying
$$\Delta B\subseteq L(\mathcal{F}_1)\setminus\{(0,0)\}=(\{0\}\times\{\pm 2,\pm 2p\})\cup ( \{1,-1\}\times \{0,\pm \sqrt{2},\pm p\sqrt{2}\}).$$
It follows that $\Delta B$ has at least $6$ and at most $8$ elements which belong to $\mathbb{Z}_3\times p\mathbb{Z}_{p^2}$. W.l.o.g, assume that $v_1=(0,0)$,
$$v_2,v_3\in (L(\mathcal{F}_1)\setminus\{(0,0)\})\cap (\mathbb{Z}_3\times p\mathbb{Z}_{p^2})=\{(0,\pm 2p),(\pm 1,0),(1,\pm p\sqrt{2}),(-1,\pm p\sqrt{2})\}$$
and
$$v_4\in L(\mathcal{F}_1)\setminus(\mathbb{Z}_3\times p\mathbb{Z}_{p^2})=\{(0,\pm 2),(1,\pm  \sqrt{2}),(-1,\pm \sqrt{2})\}.$$
Since $v_3-v_2\in (L(\mathcal{F}_1)\setminus\{(0,0)\})\cap (\mathbb{Z}_3\times p\mathbb{Z}_{p^2})$, we can assume, up to the sign, that $v_2=(1,0)$ and $v_3=(-1,p\sqrt{2})$ or $(-1,-p\sqrt{2})$. However, in both cases there is no $v_4\in L(\mathcal{F}_1)\setminus(\mathbb{Z}_3\times p\mathbb{Z}_{p^2})$ such that $v_4-v_2,v_4-v_3\in L(\mathcal{F}_1)\setminus (\mathbb{Z}_3\times p\mathbb{Z}_{p^2})$. Therefore $\mathcal{F}_1$ cannot be extended.

Similarly, assume that $\mathcal{F}_2$ could be extended by adding a new codeword $B=\{v_1,v_2,v_3,v_4,v_5\}$ on $\mathbb{Z}_5\times \mathbb{Z}_{p^2}$ satisfying
$$\Delta B\subseteq (\{0\}\times\{\pm 2,\pm 2p\})\cup (\{1,-1\}\times \{0,\pm 1,\pm p\})\cup (\{2,-2\}\times \{0,\pm \sqrt{2},\pm p\sqrt{2}\}).$$
It follows that $\Delta B$ has at least $10$ and at most $14$ elements which belong to $\mathbb{Z}_5\times p\mathbb{Z}_{p^2}$. W.l.o.g, assume that $v_1=(0,0)$,
$v_2,v_3,v_4\in (L(\mathcal{F}_2)\setminus\{(0,0)\})\cap (\mathbb{Z}_5\times p\mathbb{Z}_{p^2})$ and $v_5$ is invertible. However, no $v_2$, $v_3$ and $v_4$ exist such that $v_2,v_3,v_4,v_3-v_2,v_4-v_2,v_4-v_3\in (L(\mathcal{F}_2)\setminus\{(0,0)\})\cap (\mathbb{Z}_5\times p\mathbb{Z}_{p^2})$. Therefore $\mathcal{F}_2$ cannot be extended. \endproof

\begin{Remark}
The number of primes $p$ satisfying $p\equiv 7 \pmod{8}$ and $p<2000$ is $78$. Theorem $\ref{MaxOOC2}$ works for $59$ of them, which are listed in the following table:
\begin{center}\tabcolsep 0.005in
\begin{tabular}{|l|}
\hline $7, 23, 47, 71, 127, 151, 167, 191, 263, 271, 311, 359, 367, 383, 431, 439, 463, 479, 503, 631, 647, 719, 727,$\\ $743, 823, 839, 863, 887, 911, 919, 967, 983, 991, 1039, 1063, 1087, 1103, 1223, 1231, 1303, 1319, 1367,$\\ $1439, 1487, 1511,1543, 1559, 1567, 1583, 1607, 1663, 1759, 1783, 1823, 1831, 1847, 1871, 1951, 1999$.\\
\hline
\end{tabular}
\end{center}
\end{Remark}

\begin{Remark}\label{rek:final}
Theorem $\ref{MaxOOC2}$ gives some $(3p^2,4,1)$-OOCs and $(5p^2,5,1)$-OOCs failing to achieve the Johnson bound by missing one codeword. Actually by Theorem $\ref{prop:PS known}(2)$, a PS$(p^2)$ is known for any prime $p\equiv 3\pmod{4}$ and $3< p< 3500$, and so applying Theorem $\ref{OOC4}$ we can obtain $(3p^2,4,1)$-OOCs and $(5p^2,5,1)$-OOCs achieving the Johnson bound for such values of $p$.

How to construct a PS$(v_1v_2)$ such that $v_1\equiv v_2\equiv 3\pmod{4}$ is a long standing problem $\cite{AFGG}$. We are grateful to one of the reviewers for pointing out the reference $\cite{bp2}$. Employing similar tricks to the ones used in $\cite{bp2}$, which are more clever than what we made in the proof of Lemma $\ref{SilverPPS}$,  we found via computer search a PS$(p^2)$ satisfying the ``silver condition" for any prime $p\equiv 7 \pmod{8}$ and $p<600$ with the exception of $p\in\{71,311,367,463\}$ whose pairs are the following:
\begin{center}
\begin{tabular}{lll}
$\theta^{2i-1}\{1,\theta\}$, & $1\leq i\leq(x_1-1)/2$; \\
$\theta^{2i}\{1,\theta\}$, & $(x_1+1)/2\leq i\leq(x_2-2)/2$; \\
$\theta^{2i-1}\{1,\theta\}$, & $(x_2+2)/2\leq i\leq(x_3-1)/2$; \\
$\theta^{2i}\{1,\theta\}$, & $(x_3+1)/2\leq i\leq (p(p-1)-2)/4$; \\
$p\theta^{2i-1}\{1,\theta\}$, & $1\leq i\leq(p-3)/4$;\\
$\{\theta^{x_1},\theta^{x_2}\}$, \\
$\{\theta^{x_3},p\}$,

\end{tabular}
\end{center}
where $\theta=1+\sqrt{2}$ generates $U(p^2)/\{1,-1\}$ and
\begin{center}
\begin{tabular}{|ccccc|ccccc|}
\hline $p$ & $\sqrt{2}$ & $x_1$ & $x_2$ & $x_3$ & $p$ & $\sqrt{2}$ & $x_1$ & $x_2$ & $x_3$ \\
\hline
$7$&
$10$&
$1$ & $4$ & $17$ &

$23$&
$373$&
$1$ & $100$ & $207$ \\

$47$&
$477$&
$1$ &$208$ &$879$ &

$127$&
$1000$&
$29$ &$218$ &$7913$ \\

$151$&
$1464$&
$1$&$2626$&$7941$ &

$167$&
$20387$&
$1$&$6890$&$12075$ \\

$191$&
$17629$&
$15$&$11320$&$17061$ &

$263$&
$55909$&
$1$&$132$&$7149$ \\

$271$&
$717$&
$71$&$1286$&$1389$ &

$359$&
$98347$&
$1$&$2686$&$3713$ \\

$383$ &
$54140$ &$1$ & $16618$& $26147$ &

$431$ &
$174798$ &
$11$& $44086$ & $45461$ \\

$439$ &
$50464$ &
$3$&$7668$&$33963$ &

$479$&
$28927$&
$21$&$52840$&$106705$ \\

$503$&
$18454$&
$3$&$46438$&$67797$ &&&&& \\\hline
\end{tabular} .
\end{center}
We believe that in order to find a PS$(p^2)$ for $p\in\{71,311,367,463\}$ a few more tricks should be used. This is a possible direction for further research.

It should be pointed out that the above construction indicates that only $3$ starter pairs $\{1,\theta\}$, $\{\theta^{x_1},\theta^{x_2}\}$ and $\{\theta^{x_3},p\}$ are sufficient to produce a PS$(p^2)$, while PSs with the same parameters were constructed in $\cite{LJ2008}$ in which $(p+1)/2$ starter pairs were required.
\end{Remark}

\section{Conclusion}

This paper introduces the notion of an almost partitionable set that is a natural generalization of a partitionable set. The necessary conditions for the existence of an APS are established in Lemma \ref{nec} and we show that these necessary conditions are also sufficient for any APS of order smaller than 300 by using the Kramer-Mesner method in Theorem \ref{thm:<300}. A silver ratio construction for APS$(p,\alpha,\beta)$ with a prime $p\equiv 7 \pmod{8}$ is presented and explored based on Buratti's work in \cite{b18}. Recursive constructions are given to produce infinite families of APS$(v,\alpha,\beta)$s for any $\alpha$ and $\beta$ satisfying the necessary conditions. As applications we construct some optical orthogonal codes achieving the Johnson bound or the Johnson bound minus one.

PSs and APSs have close relationship with ZCPS-Whs, whose existence is far from being settled so far. The construction of APSs and PSs is a subject worth of deep investigation. As an equivalent object of PSs in the cyclic case, ZCPS-Wh$(v)$s with $v\equiv 1\pmod{4}$ have been investigated for more than 20 years. Unfortunately, no product theorem to deal with $v=v_1v_2$ in general is known, where $v_1 \equiv v_2 \equiv 3 \pmod{4}$, although in the special case $v_1 = v_2$ we do have direct constructions of ZCPS-Wh$(p^2)$ for prime $p\equiv 3 \pmod{4}$ and $p <3500$ (see Theorem \ref{prop:PS known}(2) and Remark \ref{rek:final}). APSs can be seen as a generalization of ZCPS-Whs. It is an interesting question whether there are more direct constructions like the silver ratio construction for APSs in Lemma \ref{Marco}.

\section*{Acknowledgments}

The authors thank the anonymous referees for their many valuable comments and suggestions that helped improve the quality of the paper. Special thanks goes to the reviewer for leading us to know the references \cite{b19,bp2}, simplifying the proofs of Lemma \ref{Nec} and Proposition \ref{prop:new}, and encouraging us to make Section 2.4 to include two infinite families of partitionable sets in abelian groups. We also gratefully acknowledge the editor Patric \"Osterg{\aa}rd for recommending Zenodo. Research of this paper was partially carried out while the second author was visiting Beijing Jiaotong University. He expresses his sincere thanks to the 111 Project of China [grant number B16002] for financial support and to the Department of Mathematics of Beijing Jiaotong University for their kind hospitality.

\end{document}